\documentclass[preprint,3p,times,twocolumn]{elsarticle}

\usepackage{graphicx,lineno}

\usepackage{amsmath,amssymb,mathtools,underscore,subfigure}
\usepackage{bbold}
\usepackage[utf8]{inputenc}
\usepackage{url}
\usepackage{varioref}
\usepackage[colorlinks]{hyperref}
\usepackage[capitalize,nameinlink]{cleveref}
\AtBeginDocument{}

\DeclareFontFamily{U}{rcjhbltx}{}
\DeclareFontShape{U}{rcjhbltx}{m}{n}{<->rcjhbltx}{}
\DeclareSymbolFont{hebrewletters}{U}{rcjhbltx}{m}{n}
\DeclareMathSymbol{\shin}{\mathord}{hebrewletters}{152}

\newcommand{\abs}[1]{{{\left | {#1} \right |}}}
\newcommand{\skipall}[1]{}
\DeclareMathOperator{\grad}{grad}
\DeclareMathOperator{\transp}{\negthinspace^\top}

\DeclareMathOperator{\sign}{sgn}

\newcommand{\eqdef}{\coloneqq}
\newcommand{\IR}{\ensuremath{\mathbb{R}}}
\newcommand{\IN}{\ensuremath{\mathbb{N}}}

\newcommand{\der}{\ensuremath{\text{d}}}

\newcommand{\dquot}[2]{\ensuremath{\frac{\der\thinspace{#1}}{\der{#2}}}}
\newcommand{\ddquot}[2]{\ensuremath{\frac{\der^2{#1}}{{#2}}}}

\DeclareMathOperator{\D}{D}
\DeclareMathOperator{\Hs}{H}

\newcommand{\ONE}{\mathbb{1}}

\renewcommand{\phi}{\varphi}
\renewcommand{\epsilon}{\varepsilon}
\newcommand{\funto}{\ensuremath \rightarrow}

\begin{document}

\begin{frontmatter}

\title{Differential Geometric Foundations for Power Flow Computations}
\journal{arXiv}

\newcommand{\authors}[1]{#1}

\authors{
\author[label1]{Franz-Erich Wolter}
\author[label1]{Benjamin Berger}
\address[label1]{Leibniz Universität Hannover}
}

\hyphenation{qua-dra-tic}
\hyphenation{span-ned}

\begin{abstract}
This paper aims to systematically and comprehensively initiate 
a foundation for using concepts from computational differential geometry 
as instruments for  power flow computing and research. At this point
we focus our discussion on the static  case, with power flow equations 
given by quadratic functions defined on voltage space with values in power 
space; both spaces have real Euclidean coordinates.  Central issue
is a differential geometric analysis of the power flow solution space
boundary (SSB) both in voltage and in power space. We present 
different methods for computing tangent vectors, tangent planes and 
normals of the SSB and the normals' derivatives. 
Using the latter we compute normal and principal curvatures.  
All this is needed  for tracing  the orthogonal  projection  of 
curves in voltage and power space onto the SSB for points on the SSB closest to given points on the curves,
thus obtaining estimates for the distance to the SSB.  
Furthermore, we present a new high precision continuation method for power flow solutions.
We also compute geodesics on the SSB or an implicitly defined submanifold thereof and, used to define
geodesic coordinates together with their Jacobians on the manifolds.
These computations might  be the most innovative and most  significant 
contribution of this paper, because this concept  provides a comprehensive 
coordinate system for submanifolds defined by implicit equations. Therefore 
while moving on geodesics described by  the geodesic  coordinates of the 
sub manifold at hand we get, via systematic  navigation guided by geodesic
coordinates, access to all feasible operation points of the system.
We propose some applications and show some properties of the Jacobian of the power flow map. 
\end{abstract}

\end{frontmatter}
%\linenumbers
\allowdisplaybreaks
\section{Introduction and Related Work}
In recent years, power transmission networks are facing increasing challenges due to increased load, decentralized energy production and the fluctuating nature of renewable energy sources. It is therefore highly relevant to develop computational methods that aid in the design, planning, and operation of such systems. Especially for operating a power network, it is important that the algorithms be fast so they can respond in real time. But also for long-term planning, this is an issue, as better algorithms allow for more finely-grained modeling and processing of larger systems. Therefore there is a continuing need for improvement of understanding and methods regarding the flow of power through a network.

Computations in the context of power flow and power grid engineering have been essentially applications of tools from numerical analysis combined with various types of classical engineering computations modified ad hoc for the equations under consideration. For those engineering problems, we have to analyze solution sets of non-linear equations, usually restricted by constraints and possibly varying with time. We are convinced of the importance of understanding the geometric structure of those solution sets for the following reasons: Theory from Riemannian and differential geometry helps with understanding the local and the global structure in a qualitative sense. This insight, along with concepts from computational geometry, yields also results of computations that are more precise and better organized. The solution sets defined by constraints are natural geometric objects, as they are Riemannian submanifolds showing geometric structures inherited from the surrounding space. Over the recent decades, computational differential geometry has developed tools for subtle and precise computations on those submanifolds \cite{wolter2011computational}; those results also have engineering implications for the problem at hand, which we intend to elaborate \cite{wolter2016MIT,wolter2017MIT,wolter2018NTU,wolter2018CGI,wolter2019NTU}. There seem to be no systematic earlier attempts to do so. 

An AC power transmission system is usually modeled as a graph. The vertices of the graph are called buses and represent generators or consumers. The edges are weighted with complex numbers and represent transmission lines with associated admittances.

The state of an alternating current power transmission system can approximately be captured by two sets of variables indexed by the set of buses, namely the complex voltage phasors $v = v^{\text{re}} + j\cdot v^{\text{im}}$ as well as the active and reactive powers $s =  p + j\cdot q$, with $p$ being active power and $q$ being reactive power. Complex variables are usually split in their real and imaginary parts in order to obtain purely real formulations. Alternatively, they may be represented by polar coordinates, but since the Cartesian representation will lead to quadratic equation systems, giving several special properties which can be exploited \cite{makarov2000properties}, we mostly prefer Cartesian coordinates. It is possible to abstract from the dynamical behavior of the system \cite{dobson1994irrelevance}, keeping only a necessary criterion for the system being in a stable state. This criterion can usually be expressed by requiring that the difference of the vector of powers and a certain function $F:\IR^n\funto\IR^n$ (called the power flow map), applied to the voltage variables, is zero. The function $F$ is quadratic in the real and imaginary parts of voltage variables and maps the voltage space into the power space. The matrices that define the quadratic form are sparse because their entries arise from transmission lines, and each bus is connected to only a few others. The power space may also contain voltage magnitude dimensions for some buses, adding the obvious quadratic terms to $F$ that compute voltage magnitudes from real and imaginary parts. Nevertheless, we will still call it power space. We also will use the letter $p$ for coordinates in power space regardless of whether the power is active or reactive, for simplification and because we abstract from the electro-technical details to arrive at a geometric formulation, where it would be confusing to use different symbols for dimensions which behave basically the same, having quadratic dependence on the voltages.
\begin{figure}[tp]
\fbox{
	\parbox{7.5cm}{
		\center
		\includegraphics[width=7cm]{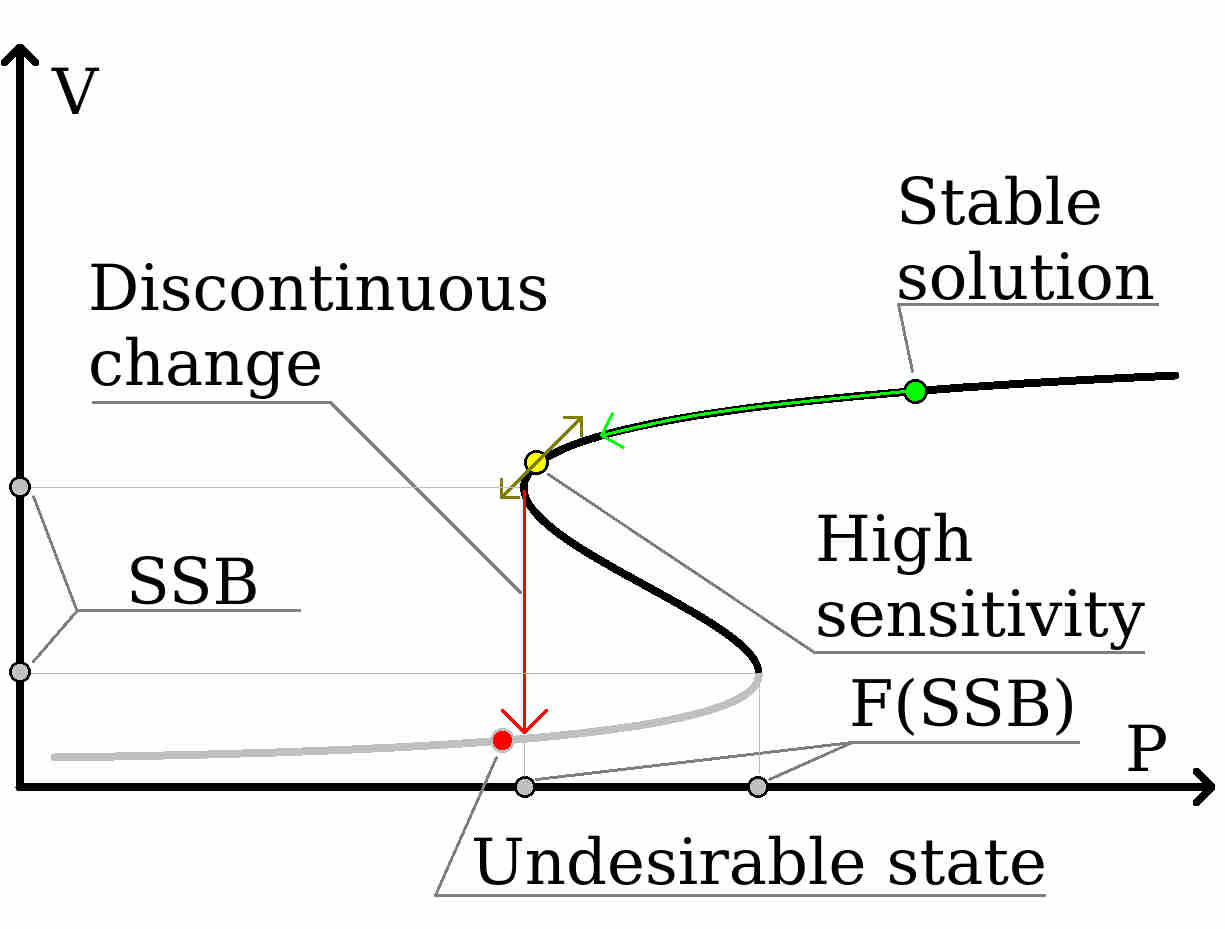}
		\caption{Schematic illustration of power flow map and SSB.}
		\label{figJump} 
	}
}
\end{figure}
The number of solutions for $v$ of the power flow equation system $p = F(v)$ varies depending on $p$. There are various algorithms for solving the power flow equations. Usually one is not interested in just any solution, but a specific one of engineering relevance characterized by high voltages and low currents.

The values of $p$ where the number of solutions changes are characterized by the existence of a $v$ with $p=F(v)$ and $\det (\D F(v))=0$, where $\D F$ denotes the differential, or Jacobian, of $F$. The set of all $p$ or of all $v$ where this is the case we call SSB (Solution Space Boundary). We will make it clear from the context whether by SSB we refer to this singular surface in voltage space or its image in power space. If $p$ changes due to varying energy generation or consumption so that it crosses the SSB, there may not exist a nearby $v$, or even none at all, that solves the power flow equations (compare \vref{figJump}). Physically, this corresponds to certain voltage collapse phenomena. The aim is then to design and operate power systems so that they are safe from collapse events and can be restored to desired states, should a contingency happen. This goal needs to be fulfilled also together with the need to minimize energy loss due to heating of transmission lines, or to minimize generation costs, or to fulfill constraints arising from market and engineering concerns.

The differential of the power flow map, $\D F$, is non-invertible if and only if its determinant is zero. The determinant is zero if and only if the smallest absolute eigenvalue of $\D F$ is zero. The differential being singular is also equivalent to the existence of a kernel vector $k$ with $(\D F) k = 0$ or a left kernel vector $(\D F)\transp\tilde k=0$.

A fairly recent overview regarding the voltage stability problem is presented in \cite{makarov2014non}. This report discusses and analyzes current methods that search for SSB points along one dimensional subspaces of power space. The report claims that the method it proposes is better than the previously available ones, but many points remain unclear, necessitating further study of the structure of the SSB (See \vref{secTopo} of this paper). The method exemplifies that a good analysis of the specific problem at hand can lead to more elegant solutions, as it exploits that the quadratic nature of $F$ implies that the Jacobian of $F$ is a linear function of the voltage, and thus can be linearly interpolated in voltage space. This is used to construct an efficient line search algorithm that is supposed to find all points along a line in voltage space where the determinant of the Jacobian vanishes, that is, points of the SSB. In turn, this enables better understanding of the SSB.

Turitsyn et al$.$ have recently analyzed efficient algorithms employing LR and QR decomposition for computing the gradients of certain matrix entries that are zero iff the determinant of the Jacobian of $F$ is zero \cite{ali2017transversality}. This provides alternative means of locating the SSB which are numerically more accurate.

Very important for our interests are previous works dealing with analysis of the SSB. Here, especially contributions of Hiskens and Dobson are noteworthy \cite{dobson1992observations,alvarado1994computation,dobson2003distance,hiskens2001exploring,hiskens1995analysis}. Dobson proposed procedures for local distance estimations in power space employing computations of principal curvatures, which is indeed a central differential geometric concept. This was a preliminary analysis having obvious limitations that we will address, hoping to propose significant improvements. Hiskens contributed a meanwhile well established continuation method for computing points and curves on the SSB. Moreover, he disproved a previously assumed convexity statement regarding the SSB in power space. Convexity properties of the domain are very important for optimization problems as well as for analysis of intersections of lines with the domain boundary. As convexity can be easily analyzed using differential geometric entities such as principal curvatures, it is certainly helpful to make systematic use of the latter to improve the aforementioned engineering computations wherever convexity questions are relevant. If convexity can be ensured for a specific case, better algorithms can be employed, and solving some problems, especially global optimization problems, may even be feasible in the convex case only.

Makarov et al$.$ showed how the quadratic nature of $F$ allows for recasting the problem of finding closest points on the SSB, and thus distances to it, as a quadratic optimization problem \cite{makarov1994continuation}.

A major part of the our intended research consists in systematically employing computational differential geometric methods in the field of electrical power engineering. This has been done by others to some degree, see for example \cite{bolognani2015fast} for an application of tangent space calculations to the power flow manifold $\{(p, v)\,\vert\, p-F(v)=0\}$ (which is much more accessible than the SSB), but not to the extent that we have in mind. Some basic facts are well-known, such as that the tangent space of the SSB in power space is spanned by the columns of the Jacobian of the power flow map $F$. From this, it follows easily that the normal vector to the SSB in power space is the left eigenvector of the Jacobian belonging to the eigenvalue $0$. There also exist various methods for computing curvatures of the SSB in power space, but not in voltage space, done by Dobson \cite{dobson1993computing,dobson1993new}. He also addresses the issue of minimal distance computations in power space. However, these computations have limitations and are incomplete, for example relying on an assumption that the SSB be not too concave in a certain sense, and finding only local minima. 
% Questions of curvatures and minimal distances are relevant to the proposed research project, especially regarding work packages 3, 5, and 6.

In the field of geometric modeling, methods of computational geometry have been in use for quite a while. Indeed, those methods need to be adapted to the case at hand, and the particular difficulties involved (High dimension, singularity, and implicit presentation of the SSB) need to be addressed. There exist classical textbooks on differential geometry presenting formulas for computing principal curvatures, normal curvatures and geodesics. The literature usually focuses on cases where manifolds are presented in an explicit way by means of a parameterization \cite{spivak1981comprehensive,manfredo1976carmo}, however, there also are books dealing with hypersurfaces defined implicitly \cite{thorpe2012elementary}. The SSB is given by such an implicit definition, but its (higher) derivatives needed for differential geometric computations are harder to obtain than with more commonly encountered implicit functions. 

An introduction to power flow related optimization problems can be found in \cite{frank2016introduction}. An overview of the algorithms employed to solve optimal power flow problems is presented in \cite{glavitsch1991optimal}. The field of Optimal Power flow (OPF) comprises several optimization problems, ranging from long-term network planning to controlling active and reactive power injections on a time scale of minutes. The goal function being minimized may for example be given by the generation costs or by the local or global transmission losses. The problem formulation includes constraints arising from engineering concerns and sometimes market requirements, in addition to the power flow equations. Numerous algorithms are employed for solving these problems. Since convex optimization problems enjoy the property that each local optimum is also a global optimum, and since the aforementioned goal functions are convex, being able to transform the problem into a convex one or to check convexity of the domain in a case under consideration is very helpful. D$.$ Molzahn has made substantial contributions in the area of optimal power flow computations, exploring conditions for when convex relaxations of the problems fail to find optimal operating points \cite{molzahn2016convex}. 

\section{Geometric Entities on the SSB}
The SSB normal $N_P$ in power space at a point is simply the left eigenvector for the eigenvalue $0$:
\begin{align}
(\D F)\transp N_P &= 0 \\
N_P\transp N_P&=1.
\end{align}
This means that the columns of $(\D F)$ are all orthogonal to the normal, so they span the tangent space.
Because $N_P$ is in the kernel of $(\D F)\transp$, we will also call it ${\tilde k}$, whereas $k$ will be the kernel of $\D F$.

In voltage space, things are more complicated. Because the SSB is an iso-surface where $\lambda_0$, the smallest eigenvalue of $(\D F)$, is zero, the normal is collinear with the gradient of $\lambda_0$. Thus the components of the not-yet-normalized normal $N_V$ are:
\begin{align}
(N_V)_i &= \dquot{\lambda_0}{V_i}.
\end{align}
The equation system 
\begin{align}
\label{ev}
(\D F - \lambda_0 \ONE)\transp  {\tilde k} &= 0 \\
{\tilde k}\transp {\tilde k}&=1, \label{knorm}
\end{align}
expresses that $(\lambda_0, {\tilde k})$ are an eigenpair. By differentiating \vref{ev}
with respect to a voltage variable $V_i$, we obtain the following linear equation system with unknowns $\dquot{\lambda_0}{V_i}$ and $\dquot{{\tilde k}}{V_i}$:
\begin{align}
\label{a}
(\D F - \lambda_0 \ONE)\transp\dquot{{\tilde k}}{V_i} - {\tilde k} \dquot{\lambda_0}{V_i}&= -\dquot{(\D F)\transp }{V_i} {\tilde k}\\
\dquot{{\tilde k}\transp}{V_i} {\tilde k}&= 0. \nonumber
\end{align}
Setting $\lambda_0=0$, but keeping its derivatives, this simplifies to
\begin{align}
(\D F)\transp\dquot{{\tilde k}}{V_i} - {\tilde k} \dquot{\lambda_0}{V_i}&= -\dquot{(\D F)\transp }{V_i} {\tilde k}\\
\dquot{{\tilde k}\transp}{V_i} {\tilde k}&= 0. \nonumber
\end{align}
Thus we not only obtain the normal vector in voltage space, but also the derivatives of the normal vector in power space (provided the rank of $\D F$ is $n-1$ and hence the zero eigenvalue is simple). Higher derivatives of the normal vector may be used be used for computing mean and normal curvatures and the shape operator \cite{do1992riemannian}, and to numerically integrate differential equations on the SSB with higher accuracy.
We are also interested in curvatures of the SSB in voltage space. So in order to get the derivatives of the unnormalized normal in voltage space, we need to differentiate \vref{a} again with respect to a voltage direction $V_j$. Then we arrive at
\begin{align}
\label{evsecder}
&(\D F - \lambda_0 \ONE)\transp\ddquot{{\tilde k}}{\der V_i\der V_j} - {\tilde k} \ddquot{\lambda_0}{\der V_i\der V_j} \\
=& -\ddquot{(\D F)\transp }{\der V_i\der V_j} {\tilde k} - \dquot{(\D F)\transp }{V_i} \dquot{{\tilde k}}{V_j} -\nonumber\\
& {\left(\dquot{\D F - \lambda_0 \ONE}{V_j}\right)}^\top\dquot{{\tilde k}}{V_i}+ \dquot{{\tilde k}}{V_j} \dquot{\lambda_0}{V_i}\nonumber\\
&\ddquot{{\tilde k}\transp}{\der V_i\der V_j} + \dquot{{\tilde k}\transp}{V_i} \dquot{{\tilde k}\transp}{V_j} = 0. \nonumber
\end{align}
By inserting the previously calculated values of $\dquot{\lambda_0}{V_i}$ and $\dquot{{\tilde k}}{V_i}$, we get a linear equation system for
the second derivatives $\ddquot{\lambda_0}{\der V_i\der V_j}$ and $\ddquot{{\tilde k}}{\der V_i\der V_j}$. Note that this equation system has the same matrix as the equation system in \vref{a}, regardless of $i$ and $j$. This means that it pays off to invert the matrix once. Note also that the matrix stays the same even if we differentiate more often; only the unknowns and the right hand side are different each time. Once the matrix has been inverted, derivatives of the same order can be computed in parallel. Note also that in our case, the second and higher derivatives of $\D F$ vanish because $F$ is quadratic, making it computationally cheaper to calculate the right hand sides than in the general case.

There is an even simpler way for obtaining the derivatives of the normal vector in power space. Observe that the columns $\dquot F {v_i}$ of $D F$ are orthogonal to the normal $N=\tilde k$:
\begin{align}
\tilde k \cdot \dquot F {v_i} =& 0. \nonumber
\intertext{Differentiating this with respect to a voltage direction $v_j$ leads to}
\dquot{\tilde k}{v_j} \cdot \dquot F {v_i} = - \tilde k \cdot \ddquot F {\partial v_i\partial v_j}.
\end{align}
This means that component in the tangent direction $(D F)(v_i)$ of the derivative of $\tilde k$ in the tangent direction $(D F)(v_j)$ can very easily be computed from the normal vector $\tilde k$ and the Hessian of $F$ (which is even constant here). These results may then be transformed linearly into a more suitable basis of the tangent space.

Having these normals and their derivatives at our disposal, we can go on to calculate normal curvatures in various tangential directions $\dot c$. Let $N$ be the normal vector in either power or voltage space. In power space, \vref{knorm} states that the normal has unit length. This simplifies the following formulas, which we give for the more general case that $N$ is not normalized.

The normal curvature $\kappa_N(\dot c)$ in the direction $\dot c$ is defined as 
\begin{align}
\kappa_N(\dot c) &= W(\dot c)\cdot \dot c, \label{wg}
\end{align}
 where $W$ is the negative differential of the unit normal, $W=- \D \frac{N}{\abs{N}}$, called Weingarten map or shape operator. 
Here we assume that $\abs{\dot c}=1$; the general case can easily be derived.
In terms of the above directional derivatives, the normal curvature is
\begin{align}
\kappa_N(\dot c)&= \left(\frac{N \left( N \cdot \dquot{N}{\dot c} \right)}{\abs{N}^3} - \frac{\dquot{N}{\dot c}}{\abs{N}}\right) \cdot \dot c.
\end{align}
If $\dot c$ is actually the tangent vector of a curve $c$, and we also know $\ddot c$, there is a simpler way to compute the normal curvature:
\begin{align}
\kappa_N(\dot c)&= \frac{N}{\abs{N}} \cdot {\ddot c}. \label{nks}
\end{align}
Only the component of ${\ddot c}$ perpendicular to the surface is relevant in \vref{nks}.

In fact, we can use this for computing the shape operator without the need do differentiate the normal vector: 
Choose $\frac{n^2-n}{2}$ tangent directions $\dot c_i$. 
Since the curves $c_i$, assumed to be parametrized proportional to arc length, run inside the surface, the normal component of each $\ddot c_i$ is uniquely determined by $\dot c$ (see \vref{ddc} for a way to compute one such $\ddot c_i$) and can be used in \vref{nks} to compute the normal curvatures in the directions $\dot c_i$. By inserting all the known normal curvatures into $\frac{n^2-n}{2}$ instances of \vref{wg}, one instance for each $i$, we get a linear equation system for the components of the matrix representation of $W$:
\begin{align}
W(\dot c_i)\cdot \dot c_i = \kappa_N(\dot c_i).
\end{align}
This is an equation system with $\frac{n^2-n}{2}$ equations and the same number of unknowns, because $W$ is self-adjoint with respect to the Riemannian metric tensor $g$, so only $\frac{n^2-n}{2}$ entries are needed to determine it. In particular $W_{ik} = -\sum_j g^{ij}L_{jk}$ where the $g^{ij}$ are the components of the inverse of the first fundamental form, and the $L_{jk}$ are the components of the second fundamental form, which form a symmetric matrix.

By cleverly choosing the directions $\dot c_i$ , we can ensure that the equation system is solvable and sparse. Not only is it sparse, but efficiently solvable by a series of substitutions.
The scheme for choosing the directions is as follows: Choose the first $n-1$ directions as the standard basis of the tangent space, and the remaining directions as the sums of each pair of distinct basis vectors.
The equation system for the components of $L$ then consists of the equations:
\begin{align} \label{Lij}
L_{ii} &= \kappa_N\left(\dot c_i\right) &\quad \text{for } 1\leq i<n \\
L_{ii} + 2 L_{ij} + L_{jj} &= \kappa_N\left(\dot c_i + \dot c_j\right) &\quad \text{for } 1\leq i< j<n   .\nonumber
\end{align}

To see how second derivatives $\ddot c_i$ of curves going in the directions $\dot c_i$ can be found, let us first consider some algebraic properties of the the Jacobian of a quadratic function $F : \IR^n \funto \IR^n$.
Such a function can always be written using symmetric matrices $A_i$, $1\leq i\leq n$:
\begin{align*}
	F(v) =  \begin{pmatrix}
		F_1(v)\\
		\vdots\\
		F_n(v)
	\end{pmatrix} = \begin{pmatrix}
	v^\top A_1 v\\
	\vdots\\
	v^\top A_n v
	\end{pmatrix},
\end{align*}
Because these matrices are symmetric, the Jacobian is simply 
\begin{align*}
\D F(v) = 2 \cdot
\begin{pmatrix}
	v^\top A_1\\
	\vdots\\
	v^\top A_n
\end{pmatrix}
\end{align*}
and we have that
\begin{align}
\label{eq27}
(\D F(x)) v &= 2 \cdot
\begin{pmatrix}
	x^\top A_1 v\\
	\vdots\\
	x^\top A_n v
\end{pmatrix}= 2 \cdot
\begin{pmatrix}
	v^\top A_1 x\\
	\vdots\\
	v^\top A_n x	
\end{pmatrix} \nonumber\\
& = (\D F(v)) x.
\end{align}
The second derivative of a quadratic form is constant, so the vector of Hessians of the components of $F$ looks like this:
\begin{align*}
\Hs F(x) = \dquot{\D F(x)}{x} = 2 \cdot
\begin{pmatrix}
	A_1\\
	\vdots\\
	A_n
\end{pmatrix}.
\end{align*}
That allows for the following transformation:
\begin{align}
\label{eq28}
v\transp\left(\Hs F(x)\right)w =& \nonumber \\
 \left(\dquot{\D F(x)}{x} v \right )\cdot w 
=& 2 \cdot
\begin{pmatrix}
	A_1 v \\
	\vdots\\
	A_n v
\end{pmatrix} \cdot w \\
=& 2 \cdot
\begin{pmatrix}
	(v\transp A_1)\transp \\
	\vdots\\
	(v\transp A_n)\transp
\end{pmatrix} \cdot w \nonumber \\
=& 2 \cdot
\begin{pmatrix}
	v\transp A_1 w \\
	\vdots\\
	v\transp A_n w
\end{pmatrix} \nonumber \\
=& (\D F(v)) w,
\intertext{thus}
v\transp\left(\Hs F(x)\right) =& (\D F)(v). \label{eq28compact}
\end{align}

With this, we can obtain the second derivatives of the curves $c_i$ in voltage space. The function F is a homogeneous quadratic form, hence the SSB is composed of cones: For any $c$ so that $\D F(c)$ is singular, the vector from the origin to $c$ is tangential to the SSB.
With $k$ being the kernel of $\D F$, we then have the relation
\begin{align}
(\D F (c)) k &= 0.
\intertext{Differentiating with respect to the curve parameter yields}
\left(\dquot{\D F (c)}{c} \dot a \right) \cdot k + (\D F (c)) \dot k &= 0, 
\intertext{then applying \vref{eq28} leads to}
(\D F (\dot c)) k + (\D F (c)) \dot k &= 0, 
\intertext{and with \vref{eq27} we get}
(\D F (k)) \dot c + (\D F (c)) \dot k &= 0.\label{l217} 
\intertext{To obtain second derivatives, we differentiate again:}
(\dquot{\D F (k)}{k} \dot k) \cdot \dot c + (\D F (k)) \ddot c +&  \nonumber\\
(\dquot{\D F (c)}{c} \dot c) \cdot \dot k + (\D F (c)) \ddot k &= 0,
\intertext{and apply \cref{eq28} and \vref{eq27} again, the former twice:}
(\D F (\dot k) ) \dot c + (\D F (k)) \ddot c + & \nonumber \\
(\D F (\dot c) ) \dot k + (\D F (c)) \ddot k &= 0.\\
(\D F (k)) \ddot c  + (\D F (c)) \ddot k &= -2 (\D F (\dot c) ) \dot k.
\end{align}
This is an under-determined equation system for the unknowns $\ddot c$ and $\ddot k$. We make it uniquely solvable by adding two conditions: First, the kernel can be required to be normalized; 
differentiating $k \cdot k = 0$ twice gives $\ddot k \cdot k = - \dot k \cdot \dot k$.
Second, we are only interested in the normal part of $\ddot c$, so we can assume without loss of generality that $\ddot c = \lambda N$ for some $\lambda\in\IR$.
This corresponds to $c$ being a geodesic.
Finally, we arrive at an equation system with $n+1$ unknowns and the same number of equations:
\begin{align}
\label{ddc}
\begin{pmatrix}
(\D F (k)) N & \D F (c) \\
0 & k\transp
\end{pmatrix}
\begin{pmatrix}
\lambda \\ \ddot k
\end{pmatrix}
=
\begin{pmatrix}
-2 (\D F (\dot c) ) \dot k \\ - \dot k \cdot \dot k
\end{pmatrix}.
\end{align}
We still need to determine $\dot k$ for the right hand side.
If we add the condition that $k$ remain normalized, expressed as $k\cdot \dot k=0$, to \cref{l217}, we can derive a linear equation system for $\dot k$:
\begin{align}
\begin{pmatrix}
\D F (c) \\
k\transp
\end{pmatrix}
\dot k=
\begin{pmatrix}
- (\D F(k))\dot c \\ 0
\end{pmatrix}.
\end{align}
These are $n+1$ equations for $n$ unknowns, but they always admit a solution since the first $n$ equations are linearly dependent, as $c$ lies on the SSB.

We now know how to compute second derivatives of geodesics, and thus normal curvatures, in voltage space. The second derivatives of the images of curves on the image of the SSB in power space are a bit more difficult to get. For this, we use the fact that the first derivative of the image $F(c)$ with respect to the curve parameter $t$ is $(\D F(c)) \dot c$. If we differentiate this again and then apply \vref{eq28}, we have
\begin{align*}
 \ddquot{F(c)}{\der t^2} &=\dquot{(\D F(c)) \dot c}{t} 
\\ &=  \left(\dquot{\D F(c)}{c} \dot c\right)\cdot \dot c + (\D F(c)) \ddot c \\
&=(\D F(\dot c))\dot c + (\D F(c)) \ddot c.
\end{align*}

When $n$ is large and we are interested in only a few normal curvatures in power space, say one in direction $(\D F(c)) \dot c$, it is more appropriate to not compute the complete shape operator using \vref{Lij}, but rather use the following formula for applying the shape operator, given in the local coordinates provided by interpreting $F$ as parameterization of power space by voltage space:
\begin{align}
W(\dot c)_i = g^{-1} L \dot c,
\end{align}
restricted to the tangent space of the SSB image in power space. 
Here $g$ is the first fundamental form related to $F$, with components $g_{ij} = \dquot{F}{V_i} \cdot \dquot{F}{V_j}$. The matrix $\tilde L$ yields the second fundamental form of the SSB-image when restricted to tangential space. 
Its components are $\tilde L_{ij}=\ddquot{F}{\der V_i \der V_j} \cdot N$.
The matrix $g$ is singular on the SSB, with the kernel vector pointing in the direction of the normal $N$. But for vectors $\dot c$ tangent to the SSB, this is irrelevant. To avoid having to invert the singular $g$, we multiply it to the left side and obtain:
\begin{align}
g (W(\dot c)) = \tilde L \dot c. \label{sswm}
\end{align}
Then, to make this linear equation system for $W(\dot c)$ more palatable to numerical solvers, we may want to ensure the right hand side does not have a normal component by multiplying it with the projection operator $(\ONE - N N\transp)$, and then add the dyad matrix $\epsilon N N\transp$ (where $\epsilon>0$) to $g$ in order to change the zero eigenvalue associated with the normal $N$ into $\epsilon$:
\begin{align}
(g + \epsilon N N\transp) (W(\dot c)) = (\ONE - N N\transp) \tilde L \dot c.
\end{align}
This also ensures that the solution vector is tangential: The right hand side is orthogonal to $N$, and the matrix has a nonzero eigenvalue belonging to the eigenvector $N$, so clearly the solution must be orthogonal to $N$, too.

Since in our context $g$ originates from the adjacency matrix $\D F$ of a mostly planar graph with node degree in practice bounded by some small $k\in\IN$, we suggest solving this equation system using an algebraic multigrid technique. This should take $O(nk^2)$ time.

We will also need principal curvatures for distance estimation purposes. Principal curvature directions are interesting later on for recognizing what directions point toward extrema of curvature that we want to avoid. If matrices of the first fundamental form $g$ and of the second fundamental form $L$ are known, we can calculate some or all of these by solving the generalized eigenvalue problem
\begin{align}
g v = \kappa_N(v) L v
\end{align}
for eigenvalues $\kappa_N(v)$ and principal directions $v$, for example using an iterative solver. Some iterative solvers can also be used if we avoid calculating the dense matrix $L$ and instead pass a procedural description of the linear map $L$ to the solver so that $L v$ can be calculated without needing an explicit representation of $L$. This procedural description may be based on \vref{sswm}.

\section{Geodesics}
In future work, we plan to use families of geodesics to parameterize submanifolds. Geodesics have the advantage that, according to the Hopf-Rinow theorem \cite{do1992riemannian,hopf1931begriff}, every point of a connected smooth manifold can be reached from every other point by a geodesic. See \cite{thielhelm2015geodesic,thielhelm} for the theory of numerical implementation of these ideas. Thus the geodesics originating from one point can be used to parameterize such manifolds. Note that the Hopf-Rinow theorem does not apply to the SSB because it is not a manifold at the origin. But for certain subsets of interest defined by algebraic equations it does apply. 

\begin{figure}[tp]
\fbox{
	\parbox{7.5cm}{
		\center
		\subfigure[The manifold]{\includegraphics[width=7cm]{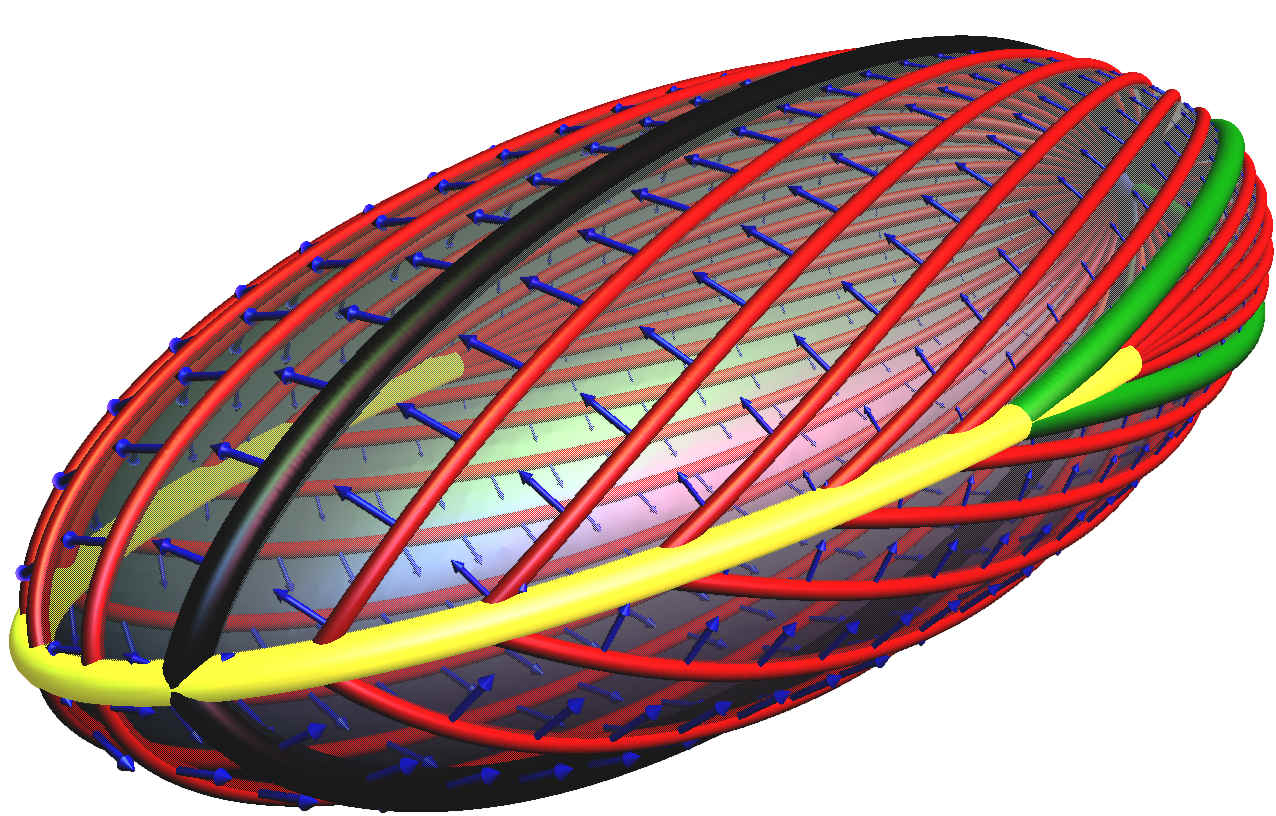}}\\
		\subfigure[The parameter space]{\includegraphics[width=6cm]{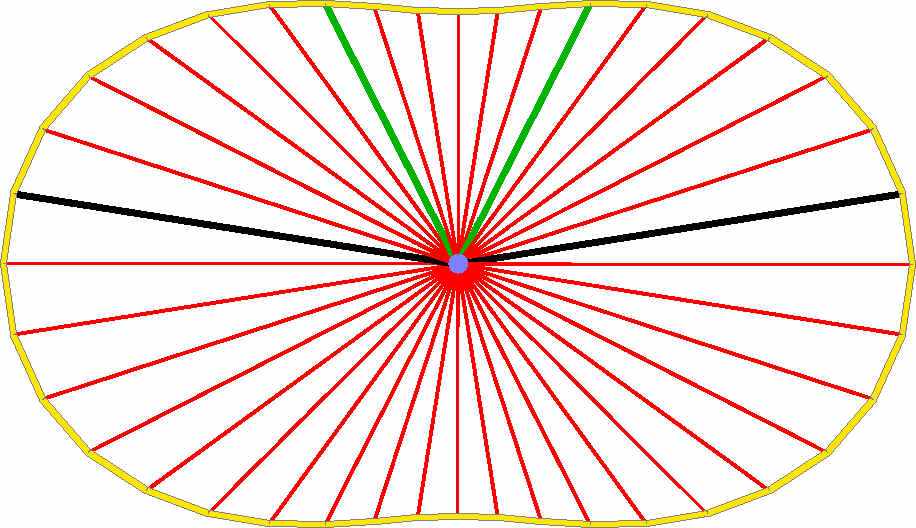}}
		\caption{A family of geodesics originating from the same point of an ellipsoid, definig a geodesic coordinate system. Geodesics are drawn in red, with two pairs of geodesics higlighted in green and black. Jacobi fields are visualized by the blue arrows and the cut locus is shown in yellow.}
		\label{figGeoEllips} 
	}
}
\end{figure}

\vref{figGeoEllips} shows a geodesic coordinate net covering the entirety of an ellipsoid along with the parameter region describing the ellipsoid in the parameter space of geodesic polar coordinates.

\begin{figure}[tp]
\fbox{
	\parbox{7.5cm}{
		\center
		\includegraphics[width=7cm]{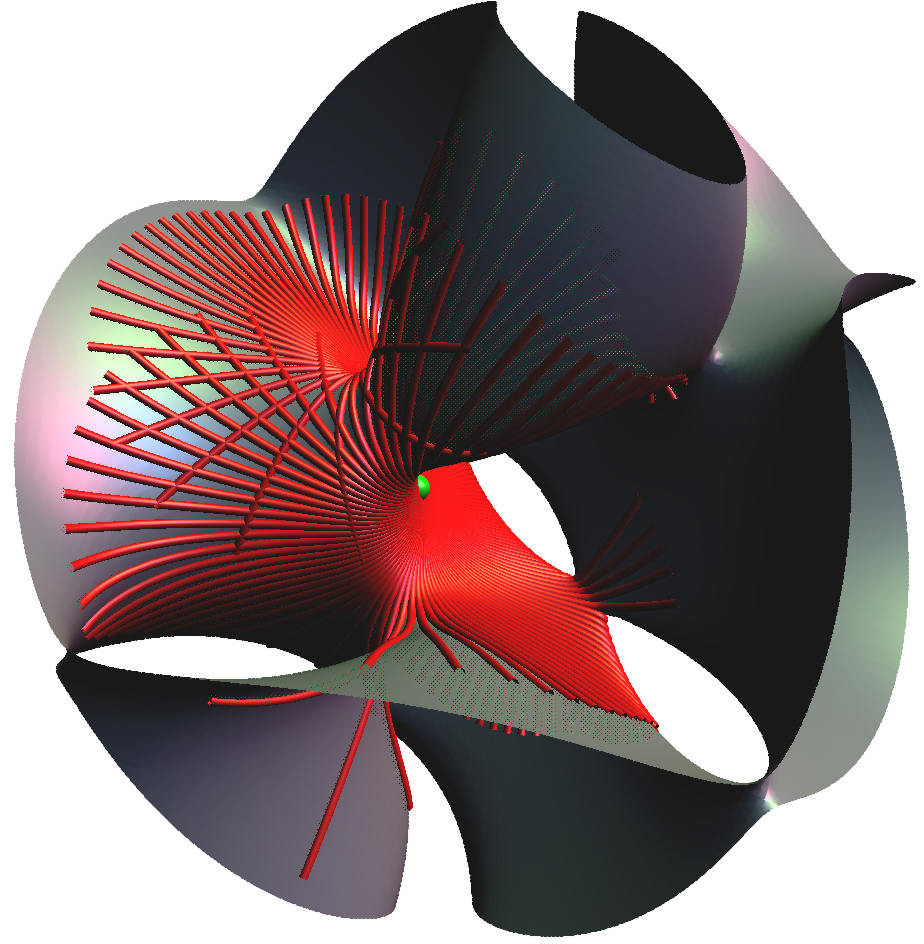}
		\caption{A set of geodesics on a 3-dimensional SSB in 4-dimensional space, visualized by orthogonally projecting everything onto the 3-sphere and then projecting stereographically into 3-dimensional space.}
		\label{figGeod4D} 
	}
}
\end{figure}

\vref{figGeod4D} shows a projection of a geodesic polar coordinate grid on an SSB for the case $n=4$. Because the SSB is composed of straight lines through the origin, it is sufficient to represent each such line by a single point. The starting directions for the geodesics have been chosen to be orthogonal to the radial direction. Note how the projections of the geodesics seem to reach every point of the surface.

\begin{figure}[tp]
\fbox{
	\parbox{7.5cm}{
		\center
		\subfigure[A geodesic coordinate net on the intersection of the unit sphere and a 3-dimensional SSB in 4-dimensional space, stereographically projected into 3-dimensional space.]{\includegraphics[width=7cm]{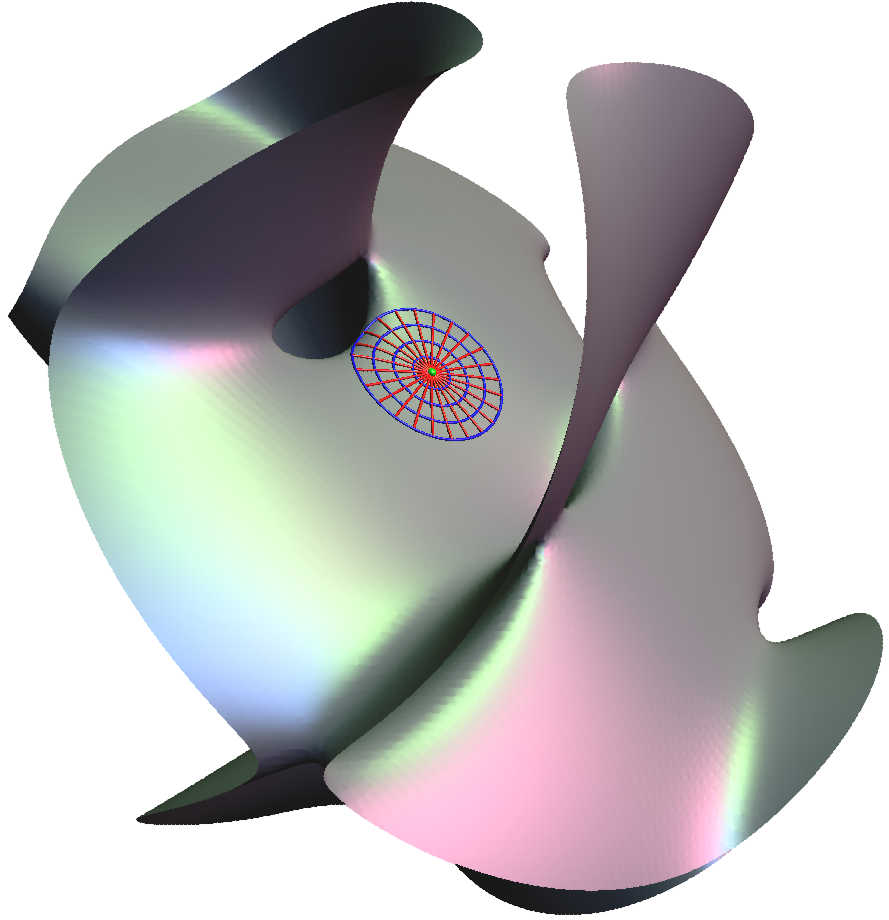}}\\
		\subfigure[The projection of the geodesic coordinate net onto the subspace spanned by coordinate directions $P_1$ and $P_2$]{\includegraphics[width=3.4cm]{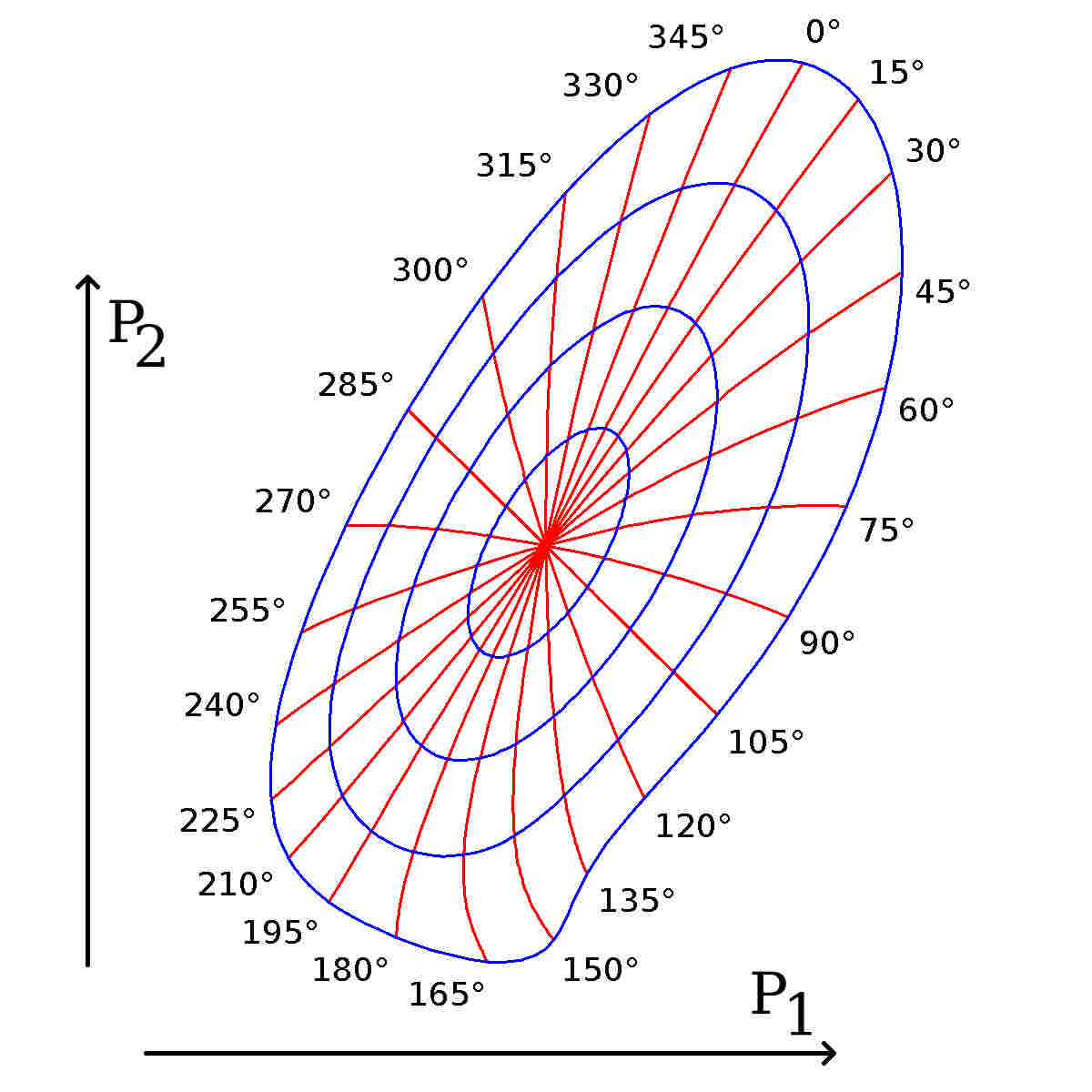}}\quad
		\subfigure[Ditto for coordinate directions $P_3$ and $P_4$]{\includegraphics[width=3.4cm]{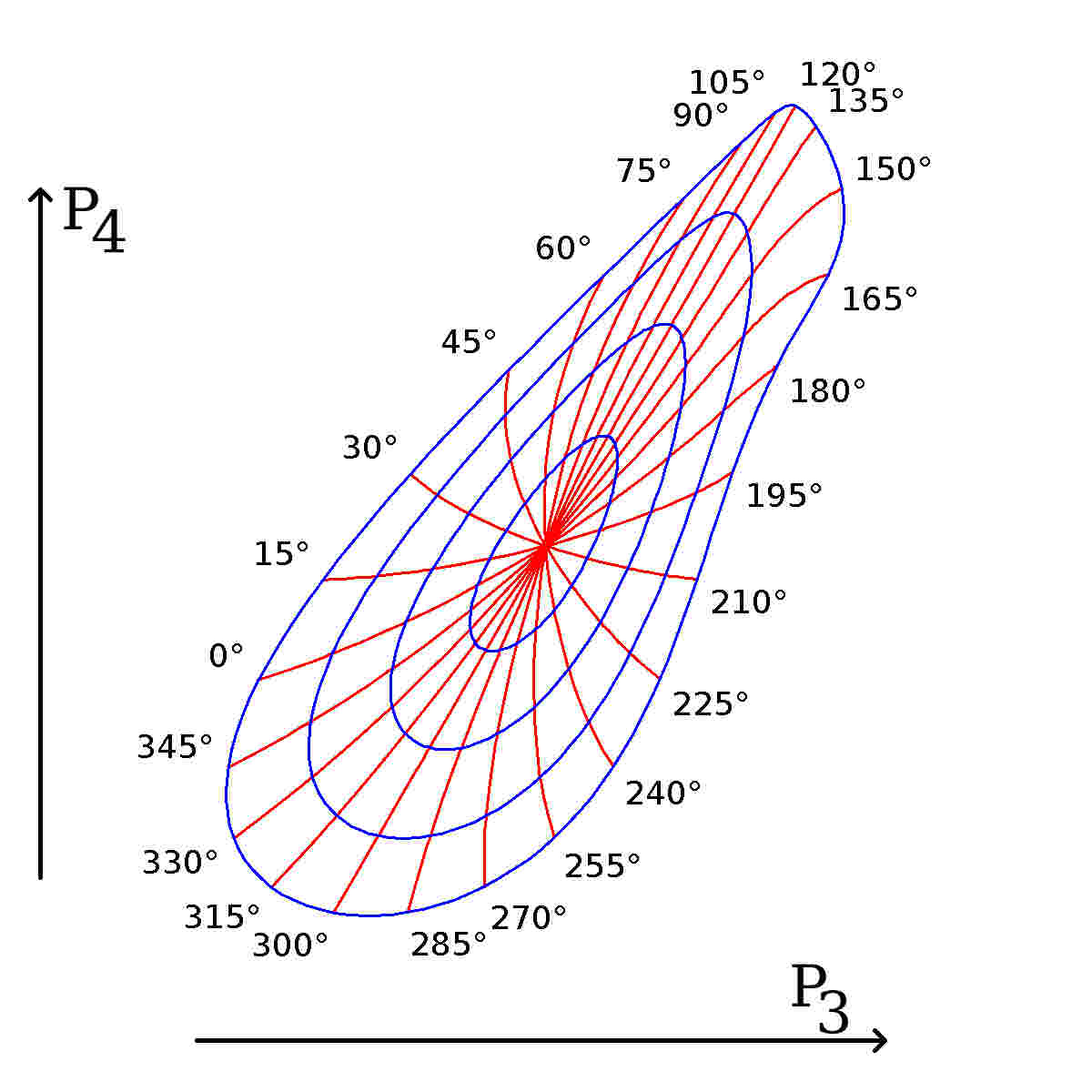}}
		\caption{}
		\label{figGeoProj} 
	}
}
\end{figure}

One factor in our motivation to study geodesics on the SSB and its submanifolds (\cite{wolter2016MIT,wolter2017MIT}, see \cref{geodHCD} on \pageref{geodHCD}) is that in certain applications, such as the question by how much some variables can be changed without the system reaching the SSB and thus becoming unstable, we are given the values of some Cartesian coordinates of a point on such an implicitly defined manifold and want to find the values of the remaining coordinates such that all coordinates together specify a point on the manifold. Our approach to this problem was first discussed in the master thesis of Gruhl \cite{gruhl}. \vref{figGeoProj} illustrates how geodesic polar coordinates may be utilized to accomplish this: If we can find the geodesic polar coordinates of a point given its $P_1$ and $P_2$ Cartesian coordinates, the coordinates $P_3$ and $P_4$ can be found by following the geodesic because the geodesic coordinate system establishes a (local) correspondence between the $(P_1,P_2)$-plane and the $(P_3,P_4)$-plane. The geodesic polar coordinates of a partially specified point can be found using an algorithm by \cite{gutschke2015differential}. 

A parameterized geodesic $\gamma: \IR \funto M$ on an implicitly defined manifold $M \subset \IR^n$ of dimension $n-1$ can be characterized by the acceleration $\ddot \gamma$ being parallel to the surface normal. If the geodesic is parameterized by arc length, this relation is given by $\ddot \gamma = \kappa_N(\dot \gamma)  N$, that is, the proportionality factor is the normal curvature in the direction of the geodesic. Here, we take $N$ to be of unit length. Since $\kappa_N$ and $N$ depend on the location $\gamma$, and $\gamma$ depends on the curve parameter $t$, this is a second order ordinary differential equation system for $\gamma$. We make it first order by replacing some instances of $\dot \gamma$ by $\alpha$: %\enlargethispage{-1.5\baselineskip}
\begin{align} \label{gjplicit}
\dot \gamma &= \alpha \\
\dot \alpha &= \kappa_N(\alpha) N. \nonumber
\end{align}
Using the start point for the initial value of $\gamma(0)$ and a unit tangential vector pointing in the initial direction of the geodesic for the initial value of $\alpha(0)$, this is an initial value problem that, when integrated for a duration of $t_1$, gives the value of the exponential map (see e.g. \cite{do1992riemannian}) at $\gamma(0)$ applied to $t_1\alpha(0)$.

%\FloatBarrier{}
\subsection{Jacobi Fields}
\begin{figure}[tp]
\fbox{
	\parbox{7.5cm}{
		\center
		\includegraphics[width=7cm]{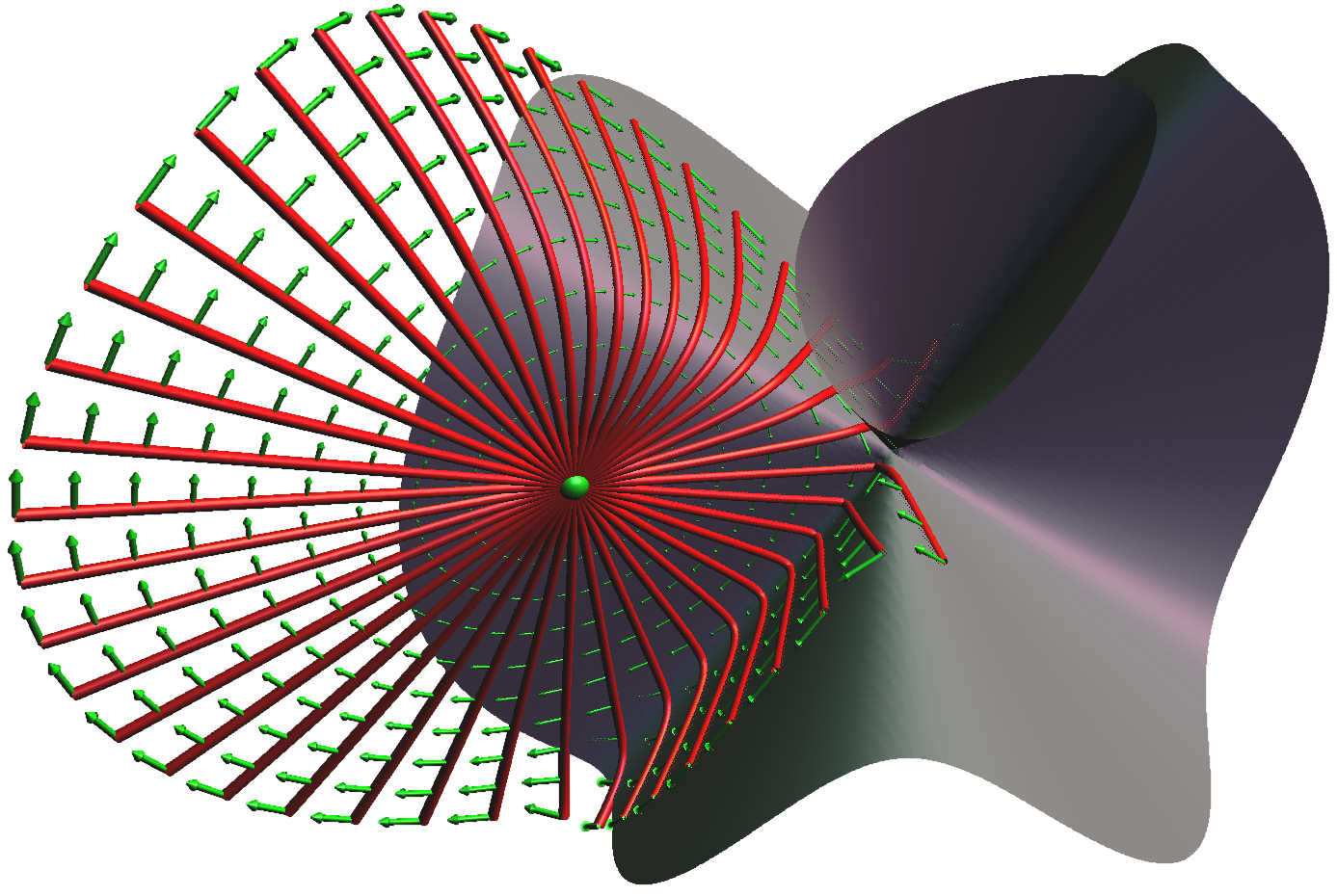}
		\caption{Illustration the exponential map and its derivative with respect to angle. Note that only part of the surface is shown, whereas the geodesics continue past the visible part.}
		\label{figGeodPolar} 
	}
}
\end{figure}

Next, we want to be able to calculate Jacobi fields; these are the derivatives of the end point $\gamma(t_1)$ with respect to a change $\alpha(t_0)^\prime$ of the starting direction $\dot \gamma(t_0) = \alpha(t_0)$. Here, $\alpha^\prime$ is the derivative of $\alpha$ with respect to some fixed variable it implicitly depends on, such as an angle. \vref{figGeodPolar} shows a geodesic polar coordinate grid (red) around a point on an SSB-like surface together with green arrows indicating the derivative of the exponential map with respect to the starting angle.

In order to avoid third derivatives of the function whose zero-set is the SSB, we use that $\kappa_N(\alpha) = \dot \alpha \cdot N$. Inserting this into \vref{gjplicit} gives us the \emph{implicit} differential equation system of \cref{Gdef} and \vref{Gdef2}. This need not worry us, because we are interested only in derivatives of its solution. What is problematic, however, is that \vref{Gdef2} really says nothing more about $\dot \alpha$ than that it should be parallel to $N$. This is the reason why later on, \vref{gjtang} is not uniquely solvable.
\begin{align} \label{Gdef}
\dot \gamma &= \alpha \\
\dot \alpha &= (\dot \alpha \cdot N) N. \label{Gdef2}
\end{align}
Differentiating this yields
\begin{align} 
\dot \gamma^\prime &= \alpha^\prime \label{gjimplicit1}\\
\dot \alpha^\prime  &= (\dot \alpha^\prime \cdot N) N \nonumber \\
	&+ (\dot \alpha \cdot (\D N) \gamma^\prime) N \qquad (=0) \label{gjimplicit2}\\
	&+ (\dot \alpha \cdot N)  (\D N) \gamma^\prime). \nonumber
\end{align}
The second term in \vref{gjimplicit2} is actually zero, as $\dot \alpha$ is parallel to $N$ (see \vref{Gdef2}), but the image of $(\D N)$ is perpendicular to $N$ because $N$ is always a unit vector. If we rearrange \vref{gjimplicit2} so that the unknown $\dot \alpha^\prime$ occurs on one side only, we get
\begin{align} 
(1-NN\transp)\dot \alpha^\prime = (\dot \alpha \cdot N) (\D N) \gamma^\prime.  \label{gjtang}
\end{align}
The matrix on the left side is the orthogonal projection operator that projects onto the tangent space. Hence the equation system is not uniquely solvable; the normal component of $\dot \alpha^\prime$ is undetermined. To fix the normal component, consider that Jacobi fields are always tangential: $\gamma^\prime \cdot N = 0$. Differentiating this twice with respect to the curve parameter yields
\begin{align}
\dot\alpha^\prime \cdot N = -2\alpha^\prime\cdot (\D N)\alpha - \gamma^\prime \cdot ((\D N)\dot \alpha + (((\D \D N)\alpha)\alpha)),
\end{align}
and so we have found the missing normal component of $\dot\alpha^\prime$. Combining this with the tangential component obtained from \vref{gjtang} gives us
\begin{align}
\dot \alpha^\prime &= (\dot \alpha \cdot N) (\D N) \gamma^\prime \\
 &- (2\alpha^\prime\cdot (\D N)\alpha + \gamma^\prime \cdot ((\D N)\dot \alpha + (((\D \D N)\alpha)\alpha))) N. \nonumber
\end{align}
Together with \vref{gjimplicit1}, this is an ODE for $\gamma^\prime$ and $\alpha^\prime$ which can be solved using standard numerical integration schemes. 

Unfortunately, we need second derivatives of the normal in order to determine the normal component of $\dot\alpha^\prime$, and that is expensive to compute although perhaps not by that much because we need only one second directional derivative ($((\D \D N)\alpha)\alpha)$ that is shared by all Jacobi fields, and we will already have computed the matrix necessary for solving \vref{evsecder} to get at the directional derivatives $(\D N) \gamma^\prime$ and can reuse that inverted matrix for finding higher derivatives. But in numerical experiments, it turned out that using the second derivative is inexact if the manifold has sharp bends, as SSBs are prone to have. To avoid relying on second derivatives, one can employ a corrector step that adjusts $\alpha^\prime$ after each step of the integrator. Observe that differentiating $\gamma^\prime \cdot N = 0$ just once with respect to the curve parameter results in
\begin{align}
\alpha^\prime \cdot N &= -\gamma^\prime\cdot (\D N)\alpha. \label{corrector}
\intertext{Additionally, we want to enforce that}
\gamma^\prime \cdot N &= 0.
\end{align}
Hence adjusting the normal components of $\alpha^\prime$ and $\gamma^\prime$ in the corrector yields (up to the expected numerical discrepancies) the same result as if we had used the correct value of the normal component of $\dot\alpha^\prime$. This variant seems to be a lot more more stable, too.

To verify our algorithm for computing geodesics, we use it to compute sectional curvatures and compare the results to sectional curvatures computed by other means.

The sectional curvature of the two-dimensional tangent subspace spanned by two orthonormal vectors $v$ and $w$ at a point $p$ is the Gauß curvature of a regular 2D surface that has $v$ and $w$ as tangents at $p$ and consists of geodesics through $p$. Hence, in an extrinsic setting, it may be computed from the Weingarten map as the determinant of the Weingarten map restricted to that subspace:
\begin{align}
\kappa_s(v, w) = \det \begin{pmatrix}
v \cdot (\D N) v & w \cdot (\D N) v \\
v \cdot (\D N) w & w \cdot (\D N) w
\end{pmatrix}.
\end{align}

The Jacobi equation provides an alternative way that utilizes a Jacobi field $J$ and its second covariant derivative along a geodesic $\gamma$ through the point in question.

The Jacobi equation \cite{do1992riemannian} states that
\begin{align}
\nabla_{\dot \gamma} \nabla_{\dot \gamma} J + R(J, \dot \gamma) \dot \gamma = 0 \label{covjac}
\end{align}
where $R$ is the Riemann curvature tensor.
The sectional curvature in the two-dimensional tangent subspace spanned by two orthonormal vectors $J$ and $\dot\gamma$ may be computed intrinsically from the Riemann curvature tensor as \cite{do1992riemannian}%TODO citation needed
\begin{align}
\kappa_s(J, \dot\gamma) = \frac{(R(J, \dot \gamma) \dot \gamma)\cdot J}{(J\cdot J)(\dot \gamma\cdot\dot \gamma) - (J \cdot \dot \gamma)^2}.\end{align}
Let us use the abbreviation $A \eqdef ((J\cdot J)(\dot \gamma\cdot\dot \gamma) - (J \cdot \dot \gamma)^2)^{-1}$. 
This together with taking the inner product of $J$ with \vref{covjac}, leads to the formula
\begin{align}
\kappa_s(J, \dot\gamma) = - A J \cdot (\nabla_{\dot \gamma} \nabla_{\dot \gamma} J).
\end{align}
The covariant derivatives can be evaluated by using that a covariant derivative is an ordinary derivative followed by a projection $P$ onto tangent space, where $P v = v - (v\cdot N) N$. Hence:
\begin{align}
\kappa_s(J, \dot\gamma) =& -A J \cdot \left (P \dquot{P \dquot{J}{t} }{t}\right).
\intertext{
The first $P$ is actually superfluous here, as the result is multiplied with a tangent vector and thus the normal part is irrelevant. So
}
\kappa_s(J, \dot\gamma) =& - A J \cdot \dquot{P \dot J }{t} \\
  =& - A J  \cdot \dquot{\dot J - \left(\dot J \cdot N\right) N}{t} \\
  =& - A J \cdot \ddot J + A \left(\dquot{\dot J \cdot N}{t}\right)(J \cdot N) \nonumber \\
  &+A\left(\dot J \cdot N\right) (J \cdot ((\D N) \dot \gamma))).
\intertext{Note that $(J \cdot N) = 0$, so:}
\kappa_s(J, \dot\gamma)  =& - A J \cdot \ddot J + A \left(\dot J \cdot N\right) (J \cdot ((\D N) \dot \gamma))) \\
=&  \frac{J \cdot \left(\left(\dot J \cdot N\right) ((\D N) \dot \gamma)) -\ddot J \right)}{(J\cdot J)(\dot \gamma\cdot\dot \gamma) - (J \cdot \dot \gamma)^2}.
\end{align}
The second derivative of the Jacobi field, $\ddot J$, is availabe from \vref{gjtang} (called there $\dot\alpha^\prime$); note that $\ddot J$ occurs only in a context where
its normal component is irrelevant, so we need not compute third derivatives in order to determine the correct normal part of $\ddot J$.

We compared both ways of computing sectional curvatures using an SSB based on a three bus system (with $n$ being $5$) by running $1000$ tests, each consisting of choosing a random point on the unit $4$-sphere, projecting it orthogonally onto the SSB and then selecting, at that point $p$, two tangential directions $v$, $w$ perpendicular to each other. The geodesic and Jacobi field ODE was integrated for the geodesic starting at $p$ in the direction $v$ and the Jacobi field corresponding to a change of the starting direction in the direction of $w$. After a geodesic distance uniformly sampled from $[0.1, 1.1]$ the integrator was halted and the sectional curvature for the tangent subspace spanned by the Jacobi field and the tangent of the geodesic at the endpoint was evaluated using both methods. 
The largest relative error (Ratio of average and difference of the two ways to compute sectional curvature) encountered was less than $2.834\cdot 10^{-9}$.

\subsection{Geodesics in Higher Codimension}
\label{geodHCD}
So far, we have only addressed geodesics on surfaces with codimension $1$, such as the SSB. We are however primarily interested in geodesics on subsets of $\IR^n$ that are defined by additional constraints, and therefore have higher codimension.

Let $C_i:\IR^n \funto \IR$ for $1\leq i \leq m < n$ be a vector of functions and let $M$ be the subset of $\IR^n$ on which all $C_i$ vanish. We assume that $M$ is regular, that is, the rank of $\D C \in \IR^{m \times n}$ is equal to $m$ on $M$. To obtain a differential equation system for a geodesic on $M$, we first make use of the fact that for every curve $\gamma$ on $M$, it holds that
\begin{align}
C(\gamma(t)) &= 0,
\intertext{and hence, by differentiating once, }
(\D C(\gamma(t))) \dot \gamma(t) &= 0
\intertext{and furthermore, by differentiating again}
\dquot{(\D C(\gamma(t))) \dot \gamma(t)}{t}  &= 0,
\intertext{that is,}
(\D C(\gamma(t))) \ddot \gamma(t) &= -\dot \gamma(t)\transp(\Hs C(\gamma(t))) \dot \gamma(t). \label{hcdNuEs}
\end{align}
This is a non-uniquely solvable linear equation system for $\ddot \gamma(t)$, because $\D C$ is not a square matrix.
Of all parameterized curves through a given point $\gamma(t)$ with a given tangent vector $\dot \gamma(t)$, a geodesic has the smallest possible value of $\abs{\ddot \gamma}$, as it by definition has no geodesic curvature and thus its curvature vector consists only of the component orthogonal to the surface. The unique solution of \vref{hcdNuEs} with minimal norm can be found by taking the pseudoinverse
of $\D C$, which is $(\D C)\transp \left((\D C)(\D C)\transp \right)^{-1}$. This leads to a second order ODE system for $\gamma$:
\begin{align}
\ddot \gamma = - (\D C)\transp \left((\D C)(\D C)\transp \right)^{-1} \dot \gamma\transp(\Hs C) \dot \gamma.
\end{align}
By introducing an extra variable set $\alpha$ for $\dot \gamma$, this can be rewritten to a first order ODE system as before.

Numerically integrating this system without a corrector step can be a bit unstable. Therefore, after each step of the integrator, $\gamma(t)$ should be adjusted so that all $C_i(\gamma(t))$ are zero, and $\alpha(t)$ should be adjusted to be of unit length and orthogonal to all $\grad C_i(\gamma(t))$.

\begin{figure}[tp]
\fbox{
	\parbox{7.5cm}{
		\center
		\includegraphics[width=7cm]{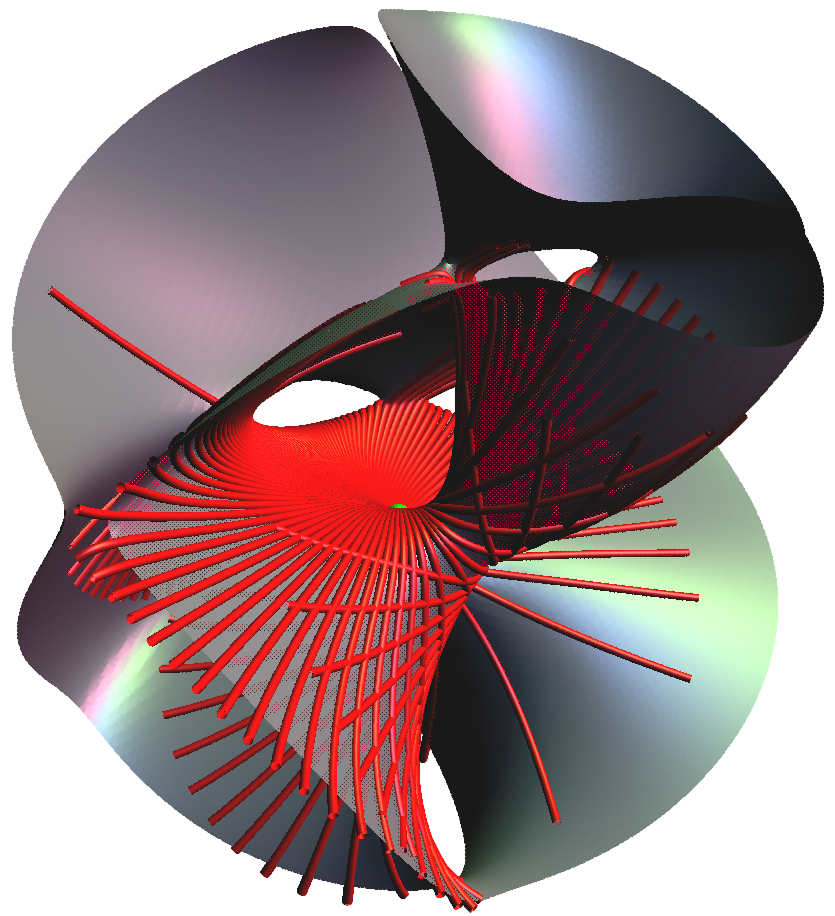}
		\caption{A set of geodesics on the intersection of the unit sphere and a 3-dimensional SSB in 4-dimensional space, visualized by projecting stereographically into 3-dimensional space.}
		\label{figGeodCd2} 
	}
}
\end{figure}

\vref{figGeodCd2} shows the stereographic projection of a $2$-dimensional surface in $\IR^4$ and some geodesics on it starting at the same point. The surface is a subset of an SSB with the additional constraint that the point has to lie on the unit sphere.

\section{Local Inversion}

Given a curve $p$ in power space, and a solution $p(0)=F(v(0))$ of the power flow equations for its starting point $p(0)$, one would often like to know the particular connected component $v$ of the preimage of $p$ in voltage space that contains $v(0)$. One application is the continuation of power flow solutions along a curve, avoiding the need to solve the equations all over again. However, standard continuation methods become unstable near the SSB, where the Jacobian is nearly singular, because they require multiplication of the inverse Jacobian by $p^\prime(t)$ to obtain $v^\prime(t)$. 

To solve this problem, we have devised two approaches \cite{hein}. Both are based on the idea of splitting the representation of the point $v(t)$ into a pair consisting of a point $q(t)$ on the SSB and a distance $d(t)$ from $q(t)$ in a particular direction $w(q(t))$, so that $v = q + d w(q)$. Depending on the approach, the direction $w(q)$ is chosen to be either the kernel of the differential at $q$, or the surface normal at $q$ (See \vref{figLI} for a sketch). The normal does not need to be of unit length, so we can simply use the gradient of the implicit function describing the SSB.

Both approaches allow us to stably compute $v^\prime(t)$ from $p^\prime(t)$ even in the vicinity of the SSB. Because the terms of $F$ are at most quadratic, its Taylor series contains at most three terms and can be used to express the dependency of $v^\prime(t)$ on $p^\prime(t)$ by a well-conditioned linear equation system. 

Both approaches should be combined: The kernel is easier to compute, needing fewer derivatives and usually providing more accuracy, but may become tangential to the SSB. If it becomes tangential, it is unreliable for the task of representing $p(t)$ and the algorithm should switch to using the normal for $w(q)$, by finding an orthogonal projection of $p(t)$ onto the SSB.

To derive the continuation method, we differentiate the equation
\begin{align}
p = F ( q + d \cdot w (q) ) 
\end{align}
with respect to $t$ (denoted by a dot over the variable):
\begin{align}
\dot p = (\D F ( q + d \cdot w(q) ) ) \left(\dot q + \dot d \cdot w(q) + d \cdot (\D w(q)) \dot q \right).
\end{align}
We want to use that $F$ is quadratic. Let $(\D F(q))$ denote the differential of $F$ at $q$ and let $(\Hs F(q))$ denote the Hessian of $F$ at $q$. Then the second order Taylor approximation is actually exact:
\begin{figure}[tp]
\fbox{
	\parbox{7.5cm}{
		\center
		\includegraphics[width=7cm]{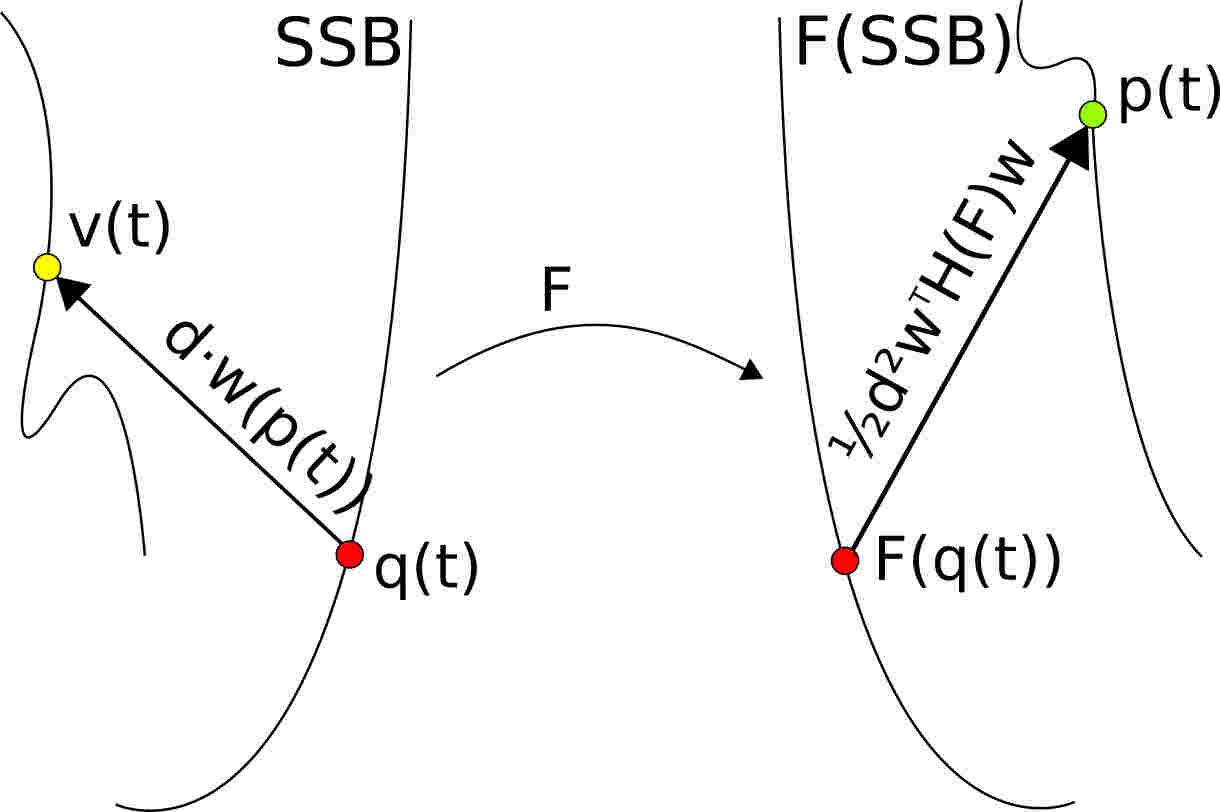}
		\caption{Calculating the local inverse of a curve using the split representation. H denotes the Hessian.}
		\label{figLI} 
	}
}
\end{figure}
\begin{align}
\label{liw}
p = F ( q ) + d \cdot (\D F(q)) w(q) + \frac{d^2}{2} \cdot  w(q) \transp (\Hs F(q))  w(q)  ). 
\end{align}
(For brevity, we will write $F$ instead of $F(q)$ and $w$ for $w(q)$.)
The linear term will vanish if we use the kernel $k(q)$ (or short $k$) of the differential at $q$ as our choice for $w(q)$, leading to the simplified equation
\begin{align}
\label{lik}
p = F + \frac{d^2}{2} \cdot k \transp (\Hs F) k.
\end{align} 
This case is illustrated in \vref{figLI}. We also provide an animation in the ancillary file {\tt{folding.mp4}} or under \url{https://www.dropbox.com/s/nlu2bp1nywmkb8c/folding.mp4} that illustrates the idea of how the power flow map acts locally like a quadratic fold at the SSB, in the generic case.
Differentiating \cref{lik} with respect to $t$ yields
\begin{align}
\label{likd}
\dot p = (\D F ) \dot q + d \dot d \cdot k \transp (\Hs F) k + d^2 \cdot k \transp (\Hs F) (\D k) \dot q.
\end{align} 
Note that we used that the derivatives of the Hessian vanish in our case. Note also that here $(\D k) \dot q$ is simply notation for the directional derivative of $k$ in the tangential direction $\dot q$. It cannot actually be computed because $k$ cannot be differentiated in all directions, as the derivative it does not even exist for nontangential directions. Later on, to solve \vref{trassDumm}, we will need to be able to apply $\D k$ to nontangential vectors. Thus, the vector field $k$ needs to be locally extended away from its natural domain, the SSB. This can be done for example by making it constant along some arbitrary transversal direction or by taking $k$ to be the right eigenvector belonging to the smallest absolute eigenvalue. A more elegant workaround is presented with \vref{trassSchlau}.

We can turn \vref{likd} into a linear equation system for $\dot d$ and $\dot q$:
\begin{align}
\begin{pmatrix}
(\D F)+ d^2 \cdot k \transp (\Hs F) (\D k) & d \cdot k \transp (\Hs F) k
%(\D F)+ d^2 \cdot w \transp (\Hs F) (\D w) & d \cdot w \transp (\Hs F) w
\end{pmatrix}
\begin{pmatrix}
\dot q \\ \dot d
\end{pmatrix}
=
\dot p.
\end{align} 
This equation system contains $n+1$ unknowns, but consists of only $n$ equations. To make it uniquely solvable one may add a single equation expressing that $\dot q$ be tangential. Then the equation system reads:
\begin{align}
\begin{pmatrix}
(\D F)+ d^2 \cdot k \transp (\Hs F) (\D k) & d \cdot k \transp (\Hs F) k \\
%(\D F)+ d^2 \cdot w \transp (\Hs F) (\D w) & d \cdot w \transp (\Hs F) w \\
N \transp (q) & 0
\end{pmatrix}
\begin{pmatrix}
\dot q \\ \dot d
\end{pmatrix}
=
\begin{pmatrix}
\dot p \\ 0
\end{pmatrix}. \label{trassDumm}
\end{align} 
This however requires the normal $N$ to be computed. We can avoid this cost by instead adding $n$ variables, namely the components of $\dot k$, and $n+1$ equations expressing that $k$ remains the kernel and of constant length as it changes along the curve by an infinitesimal amount $\dot k$. We obtain these equations by differentiating the condition that $k$ be a kernel vector of a given length with respect to the curve parameter, which is only possible if $q$ moves tangentially along the SSB, leading to the equation system
\begin{align}
\begin{pmatrix}
\D F+ d^2 k \transp (\Hs F) (\D k) & 0 & d k \transp (\Hs F) k \\
k\transp (\Hs F) & \D F & 0 \\
0 & k\transp & 0
\end{pmatrix}
\begin{pmatrix}
\dot q \\ \dot k \\ \dot d
\end{pmatrix}
=
\begin{pmatrix}
\dot p \\ 0 \\ 0
\end{pmatrix}.
\end{align}
Now that we have $\dot k$ as an explicit variable, we can also use that $\dot k = (\D k) \dot q$ to get rid of the ill-behaved $\D k$. The equation system now is
\begin{align}
\begin{pmatrix}
\D F& d^2 \cdot k \transp (\Hs F) & d \cdot k \transp (\Hs F) k \\
k\transp (\Hs F) & \D F & 0 \\
0 & k\transp & 0
\end{pmatrix}
\begin{pmatrix}
\dot q \\ \dot k \\ \dot d
\end{pmatrix}
=
\begin{pmatrix}
\dot p \\ 0 \\ 0
\end{pmatrix}. \label{trassSchlau}
\end{align}

Note that in our case $F$ is a quadratic function. According to \cref{eq28compact}, $k\transp (\Hs F)$ evaluated anywhere equals $\D F$ evaluated at $k$. This simplifies the formulas used here in several places.

We use this equation system to implement a continuation method that traces how the preimage $v = q + d \cdot k$ of $p$ evolves as $p$ changes in the direction $\dot p$. There are some problems with this that need to be addressed: First, as the curve $q$ is constructed, it may deviate from the SSB due to numerical inaccuracies. If that happens, it should be corrected by projecting it back onto the SSB along the direction $w(q)$. Second, we found in experiments that the simplification employed in \cref{lik}, while valid when using exact arithmetic, leads to numerical errors that can be dramatically reduced in exchange for the rather small effort of using the full \cref{liw}. Third, it may not always be possible or reliable to represent $v$ as $q + d \cdot k(q)$, namely if the kernel is (nearly) tangential. In that case, it would be more appropriate to use the normal instead of the kernel for $w$. A suitable point $q$ so that $v$ can be expressed as $q + d \cdot N(q)$ can be found using orthogonal projection (see \vref{secOrtho}). Using the normal should be avoided when the kernel is available for representing $v$ because the differential of the normal is much more expensive to compute.

For these reasons, we should also derive the equation system for $\dot q$ and $\dot d$ when the full \vref{liw} is used. Differentiating \cref{liw} with respect to $t$, we get:
\begin{align}
\label{liwd}
\dot p =& (\D F) \dot q \\
+& \dot d \cdot (\D F) w + d \cdot \dot q \transp (\Hs F) w\nonumber \\
& + d \cdot (\D F) (\D w) \dot q \nonumber \\
+& d \cdot w \transp (\Hs F) k + d^2 \cdot w \transp (\Hs F) (\D w) \dot q \nonumber.
\end{align} 
Proceeding as above, we can again derive a linear equation system from this that tells us $\dot q$ and $\dot d$, given $\dot p$

If first order derivatives should happen to provide insufficient accuracy for numerical integration, second order derivatives can be employed. 
%TODO How?
 Differentiating \vref{liwd} again yields
\begin{align}\label{liwdd}
\ddot p =& \dot q\transp (\Hs F) \dot q + (\D F) \ddot q \\
+& \ddot d \cdot (\D F) w + \dot d \dot q \transp \cdot (\D F) w + \dot d \cdot (\D F) (\D w) \dot q \nonumber \\
& + \dot d \cdot \dot q \transp (\Hs F) w + d \cdot \ddot q \transp (\Hs F) w+ d \cdot \dot q \transp (\Hs F) (\D w) \dot q \nonumber \\
& + \dot d \cdot (\D F) (\D w) \dot q + d \cdot (\dot q\transp (\Hs F)) (\D w) \dot q \nonumber \\
& \quad + d \cdot (\D F) (\dot q \transp\D w) \dot q + d \cdot (\D F) (\D w) \ddot q \nonumber \\
+& \dot d \cdot w \transp (\Hs F) w + d \cdot w \transp (\Hs F) (\D w) \dot q + d \cdot w \transp (\Hs F) w \nonumber \\
& \quad + 2 d \dot d \cdot w \transp (\Hs F) (\D w) \dot q + d^2 \cdot ((\D w) \dot q) \transp (\Hs F) (\D w) \dot q  \nonumber \\
& \quad + d^2 \cdot w \transp (\Hs F) (\dot q \transp (\Hs w) \dot q) + d^2 \cdot w \transp (\Hs F) (\D w) \ddot q \nonumber .
\end{align} 
The condition that $\dot q$ be tangential gives us the linear constraint $ N(q) \cdot \ddot q = -((\D N(q))\dot q) \cdot \dot q$, which together with \vref{liwdd} forms a linear equation system for $\ddot q$ and $\ddot d$.

Note that \vref{liwdd} contains the term $\dot q \transp (\Hs w) \dot q$. The Hessian of $w$ can be very expensive to compute, especially if $w=N$. Again, however, this is just notation for the second directional derivative of $w$ in the direction $\dot q$. This can be found according to \vref{evsecder} by taking the derivatives with respect to $V_i$ and $V_j$ both to be directional derivatives for direction $\dot q$.

In the special case that we are inverting a curve that runs entirely within the SSB, we can specialize equation system \vref{trassSchlau} by setting $d=0$ and removing the unknown $\dot d$ along with the third column of the matrix. To get at the second derivatives in that special case, we can then also differentiate this system with respect to $t$ to get an equation system for $\ddot q$ and $\ddot k$ in terms of $F$, $\ddot p$, $q$, $\dot q$, $k$ and $\dot k$. It turns out that this equation system uses the same matrix as the one it was derived from.

\subsection{Numerical Results}
We tested the local inversion algorithm on power network configurations from the power system test case archive of the University of Washington \cite{tcarchive}. In particular, we compared the precision of three of our methods: The method that uses the kernel and omits the linear term of the Taylor expansion because it is theoretically zero, the method that uses the kernel but includes the linear term, and the method using the normal vector. We start with curves in voltage space, map them through $F$ into power space and invert the results back into voltage space using the three algorithms. By comparing the result with the original curve, we can estimate the precision, which is shown in \vref{tabPrecision} taken from the defence presentation of the master thesis of Hein \cite{hein}.

\begin{table*}[hbt]
	\centering
	Step size $10^{-9}$ \\
	\begin{tabular}{|c|c|c|c|}
		\hline 
		\textbf{\# of buses} & \textbf{Kernel w/o lin.} & \textbf{Kernel with lin.} & \textbf{Normal}\\ 
		\hline 
		14 & $2.38774883365\cdot 10^{-5}$ & $2.37706864821\cdot 10^{-5}$ & $2.14982954847\cdot 10^{-5}$\\ 
		\hline 
		30 & $7.05835279589\cdot 10^{-9}$ & $7.05834973039\cdot 10^{-9}$ & $7.97487174331\cdot 10^{-7}$ \\ 
		\hline 
		57 & $1.64770875908\cdot 10^{-7}$ & $1.64770753583\cdot 10^{-7}$ &$3.46424134738\cdot 10^{-6}$ \\ 
		\hline 
		118 & $6.71182537764\cdot 10^{-8}$ & $1.07999417455\cdot 10^{-13}$ & $1.9659196324\cdot 10^{-7}$\\
		\hline 
	\end{tabular} \\
	Step size $10^{-7}$ \\
	\begin{tabular}{|c|c|c|c|}
		\hline 
		\textbf{\# of buses} & \textbf{Kernel w/o linear} & \textbf{Kernel with linear} & \textbf{Normal}\\ 
		\hline 
		14 & $2.41940697575\cdot 10^{-5}$ & $2.38882856096\cdot 10^{-5}$ & $1.20966340549 \cdot 10^{-4}$\\ 
		\hline 
		30 & $7.05835852083 \cdot 10^{-6}$ & $7.05835718598\cdot 10^{-6}$ & $7.97486759134 \cdot 10^{-5}$ \\ 
		\hline 
		57 & $1.64760675278 \cdot 10^{-4}$ & $1.6477270041\cdot 10^{-5}$ & $3.34938440007 \cdot 10^{-4}$ \\ 
		\hline 
		118 & $5.37431102992 \cdot 10^{-7}$ & $1.1515965505\cdot 10^{-10}$ & $2.0439378252 \cdot 10^{-4}$\\
		\hline 
	\end{tabular}
%	\caption{Durchschnittliche Distanz der invertierten Kurve zur Ausgangskurve für 14-,30-,57- und 118-Bussystemen nach $100$ Schritten mit einer Schrittweite von $10^{-9}$ (oben) und $10^{-7}$ (unten). }
	\caption{Average distance between original curve and inverted curve after 100 steps with step size $10^{-9}$ (upper table) and $10^{-7}$ (lower table) }
	\label{tabPrecision}
\end{table*}

Using the linear term in the kernel method pays off in precision, for almost no additional runtime cost.

A comparison of execution times for a single step is displayed in \vref{tabRuntime}. 
\begin{table}[hbt]	
	\centering
	\begin{tabular}{|c|c|c|}
		\hline 
		\textbf{\# of buses} & \textbf{Kernel} & \textbf{Normal}\\ 
		\hline 
		14 & 0.278 & 0.324\\ 
		\hline 
		30 & 3.329 & 6.492 \\ 
		\hline 
		57 & 42.047 & 81.992\\ 
		\hline 
		118 & 225.323 & 439.380\\
		\hline 
	\end{tabular}
	\caption{Runtime for a single step of the algorithm in seconds on a Laptop with an Intel i7 7700HQ processor }
	\label{tabRuntime}
\end{table}
Since the normal method is slower and on average less precise, the kernel method is preferable as long as the angle between kernel and tangent space is large enough.

\section{Orthogonal Projection}
\label{secOrtho}
\begin{figure}[tp]
\fbox{
	\parbox{7.5cm}{
		\center
		\includegraphics[width=7cm]{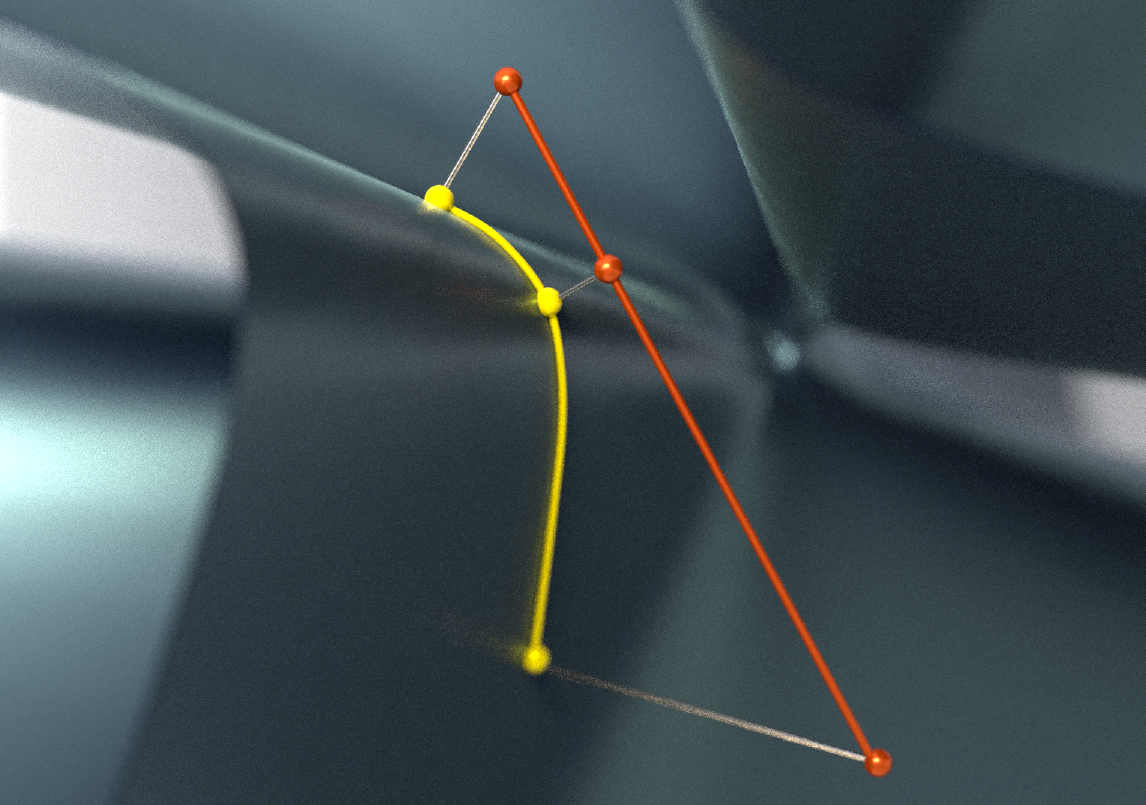}
		\caption{Orthogonal projection of a curve onto an SSB.}
		\label{figOproCurve} 
	}
}
\end{figure}

For some purposes, such as our local inversion technique when using the normal vector, it is required that, given a point $v$, we find a point $q$ on the SSB so that the normal at $q$ points along the direction $v-q$. Since the SSB is the implicit surface where  the absolutely smallest eigenvalue $\lambda$ of the power flow differential $\D F$ is $0$, we can use the approach described in the following to find that orthogonal projection $q$. Also of interest are orthogonal projections of entire curves onto the SSB \cite{pegna1996surface} (see \vref{figOproCurve}).

\subsection{An Algorithm for Computing Orthogonal Projections of Points onto the SSB}
\begin{figure}[tp]
\fbox{
	\parbox{7.5cm}{
		\center
		\includegraphics[width=7cm]{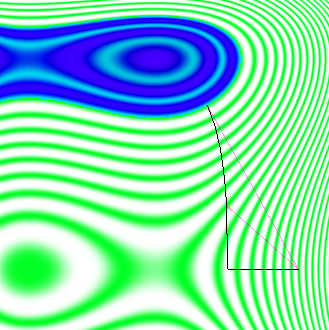}
		\caption{2D-illustration of the ODE for finding orthogonal projections onto a line (that is not SSB-like). The isolines have been visualized.
		The normal to the isoline at $t=0$ as well as a normal to an isoline halfway between $t=s$ and $t=0$ are marked.
		The black line is the curve $r$.}
		\label{figOpro2d} 
	}
}
\end{figure}
Let $\lambda:\IR^n\funto \IR$ be a function whose zero set is the surface that we want to project onto, such as the SSB with $\lambda$ being the smallest eigenvalue of $\D F$, and let $v$ be the point that we want to project.
First, we evaluate $s = \lambda(v)$. Then, we compute a parameterized curve $r$ so that $r(s) = v$, $r(0)=q$ and for all $t$ between $0$ and $s$, the normal to the iso-surface with $\lambda(r(t))=t$ and the direction $r(t) - v$ are linearly dependent. See \cref{figOpro2d} for an illustration of the idea.
The condition for lying on the isosurface is equivalent to the existence of a kernel $k$ of $\D F - \lambda \ONE$. 
We make the generic assumption that the gradient of $\lambda$ does not vanish, so that $r(t)-v$ can always be written as a multiple of it:
\begin{align} 
d(t) \cdot (\D \lambda(r(t)))\transp & = v - r(t) .
\end{align}
This is always valid at $t=s$, with $d(t)=0$. 
The expressions $\D \lambda$ (and its directional derivatives needed below) may be calculated analogous to \vref{a}.
Then, we trace the curve $r$ from $t=s$ to $t=0$ while maintaining that this relation stays true. For this, we differentiate with respect to $t$ (denoted by a dot on top of symbols). Notating the dependencies on $t$ is omitted:
\begin{align} 
\dot d \cdot (\D \lambda(r))\transp + d \cdot (\Hs \lambda(r)) \dot r &= -\dot r .
\end{align}
This is a linear equation system for $\dot r$ and $\dot d$, but it has one variable more than it has equations. The missing equation needed for unique solvability is the condition that the derivative of $\lambda$ with respect to $t$ equal $(\sign{s})$, by construction. Inserting this condition, we arrive at:
\begin{align} \label{orthoODE}
\begin{pmatrix}
d \cdot (\Hs \lambda(r)) + \ONE & (\D \lambda(r))\transp \\
(\D \lambda(r)) & 0
\end{pmatrix}
\begin{pmatrix}
\dot r \\ \dot d
\end{pmatrix}
= \begin{pmatrix}
0 \\ \vdots \\ 0 \\ \sign{s}
\end{pmatrix}.
\end{align}
We integrate this ordinary differential equation system backwards (or forwards, depending on the sign of $s$) from $t=s$ to $t=0$ and thus obtain the orthogonal projection $q=r(0)$ onto the $\lambda=0$ surface.

There is, however, one problem with this way to find the orthogonal projection of a point. In the case that $v$ lies beyond a focal surface of the SSB, there will be a $t$ so that $v$ lies \emph{on} the focal curve of the iso-surface with $\lambda(r(t))=t\cdot \sign s$, as she focal surface changes continuously with $t$ and at $t=0$ the point $v$ lies beyond it whereas at $t=s$ it does not. The numerical consequence is that the derivative of $r$, as given be \cref{orthoODE}, may diverge. If the ODE does not adapt its step size, this leads to inaccurate results, whereas if it does adapt, the step size may appoach zero and the solver gets stuck. The solution is to limit the step size from below and use \cref{orthoODE} as the predictor in a kind of predictor-corrector algorithm. We found that a good choice for the corrector step is to replace $r(t)$ with $\tilde r(t)$, which is found by minimizing $\abs{\tilde r(t)-v}^2$ under the constraint $\lambda(\tilde r(t))=t\cdot \sign s$ using projected gradient descent initialized with $r(t)$. Afterwards, $d$ needs to be updated as well, according to $d \eqdef \tilde d = \frac{N\cdot(v - \tilde r (t))}{N \cdot N}$. With a high-order adaptive step size control scheme (We used Dormand-Prince), the corrector usually has nothing to do because the predicted point already is very accurate. Only when the step size limitation becomes relevant does it have to do a few iterations. The corrector step does not use the Hessian, so it should be much cheaper than the predictor.

\begin{figure}[tp]
\fbox{
	\parbox{7.5cm}{
		\center
		\includegraphics[width=7cm]{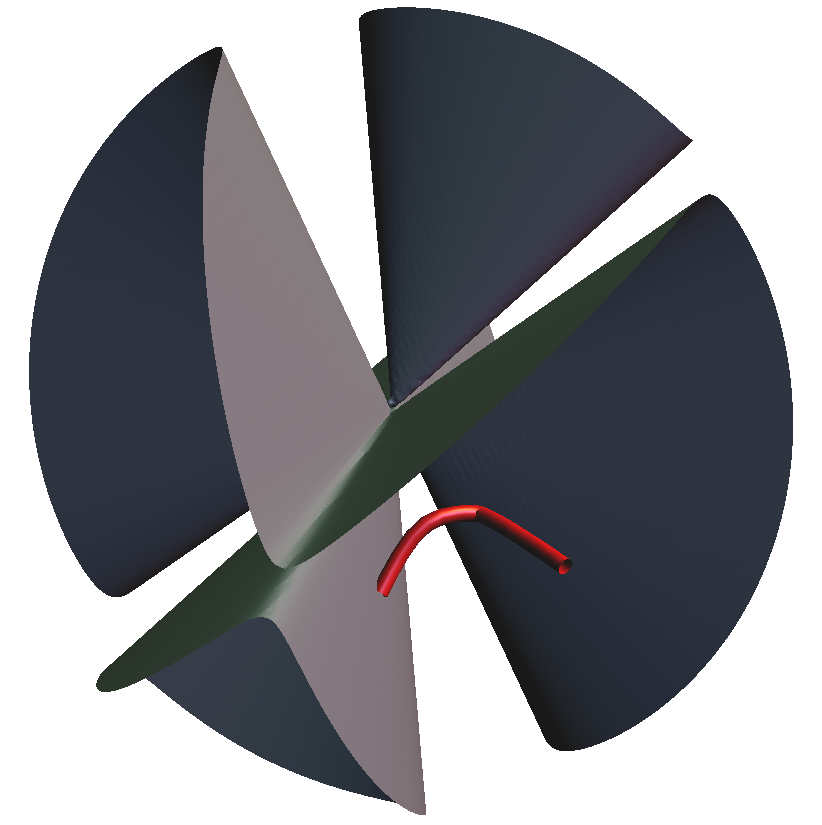}
		\caption{3D-example of the ODE for finding orthogonal projections onto an SSB-like surface.}
		\label{figOpro3d} 
	}
}
\end{figure}

\vref{figOpro3d} shows an example run of the algorithm for a surface of the form $\det \D F=0$ for a quadratic function $F:\IR^3\funto\IR^3$. It is the same function as in \vref{figDetNondegen0}.
%Unfortunately, numerical experiments indicate that that the method does not find the globally closest point.
\subsection{Tracing the Orthogonal Projection of a Curve}
Once we have found the orthogonal projection of a point $v$, we can compute its derivatives with respect to changes in $v$, and hence we can trace the orthogonal projection of a differentiable curve. So, let $v$ and $q$ now again be curves parameterized by curve parameter $t$ (unrelated to the use of $t$ in the previous paragraph), and let $d$ be a real valued function of $t$ so that $v = q + d N$, where $N(q) = \frac{\grad \lambda(q)}{\abs{\grad \lambda(q)}}$ this time is the the \emph{unit} normal to the surface with $\lambda(q)=0$. Then $d = (v - q)\transp N$, and by differentiating 
\begin{align} 
d(t) \cdot N(q(t)) & = v(t)-q(t)
\end{align}
with respect to $t$, we get
\begin{align} 
\label{orthoDotd}
\dot d  \cdot N(q) + d \cdot (\D N) \dot q & = \dot v - \dot q 
\end{align}
where 
\begin{align} 
\dot d & =   (\dot v - \dot q)\transp N(q) +  (v - q)\transp (\D N(q)) \dot q .
\end{align}
Here, some terms vanish, as $\dot q$ is orthogonal to $N(q)$ and $(v - q)$ is parallel to $N(q)$ and thus orthogonal to the image of $\D N(q)$, which is the tangent plane:
\begin{align} 
\dot d & =   \dot v \transp N(q).
\end{align}
Inserting this into \vref{orthoDotd}, we get a linear equation system for $\dot q$, given $\dot v$:
\begin{align} 
(d \cdot (\D N(q)) + \ONE) \dot q & = \dot v  - (\dot v \transp N(q)) N(q).
\end{align}
\skipall{
Note that the first term in the matrix, $N A$, is a dyad and the others are sparse. Numerical algorithms might benefit from this.
We can use this equation system to trace $q$ as $v$ changes by employing standard continuation techniques.
%TODO where does this belong? Not here.
}

\begin{figure}[tp]
\fbox{
	\parbox{7.5cm}{
		\center
		\includegraphics[width=7cm]{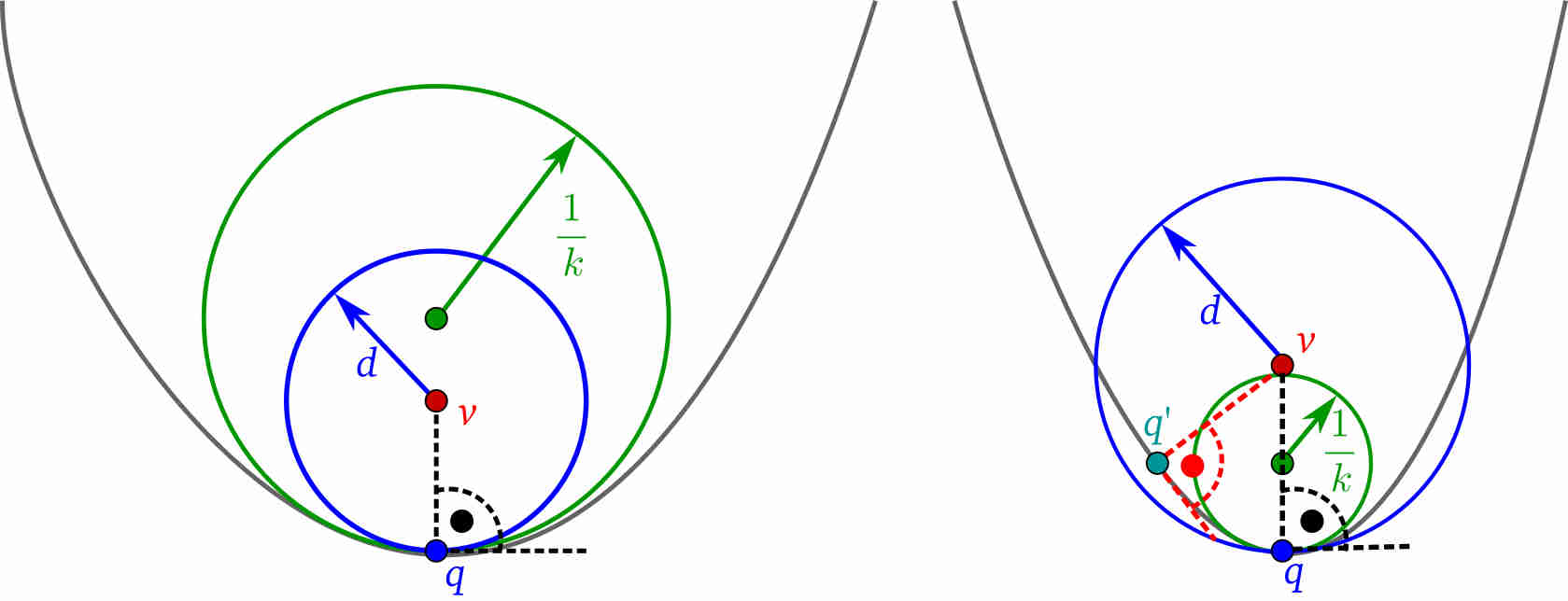}
		\caption{Illustration of the necessary condition for an orthogonal projection $q$ of $v$ to be a closest point: If the distance $d$ is larger than a radius of curvature $\frac{1}{k}$ (shown on the right), there is always another point that is closer to $v$ than $q$ is, for example the point $q^\prime$. Here, $q^\prime$ is also an orthogonal projection and a globally closest point, but already points infinitesimally adjacent to $q$ are closer to $v$ than $q$ is, because in second order they follow the osculating green circle of radius $\frac{1}{k}$, each of whose points except $q$ lies inside the blue circle of radius $d$ around $v$.
		 }
		 \label{figKruku}
	}
}
\end{figure}

The function $d$ tells the distance between $v$ and $q$. Provided that it is smaller than the smallest positive\footnote{``Positive'' means curving towards $v$; We explain here only the case where $d\leq 0$ and therefore the normal points away from $v$.} radius of curvature of the SSB at $q$, $q$ is a point on the SSB locally closest to $v$. The smallest positive radius of curvature is the reciprocal of largest eigenvalue of the Weingarten map. 
If $d$ is larger than the smallest radius of curvature, then $q$ cannot possibly be the closest point to $v$ on the SSB. See \vref{figKruku} for an illustration, adapted from the master thesis of Gruhl \cite{gruhl}.

\skipall{
%TODO Was ist hiermit?
\section{Minimum Distance Computations}
Computing the global minimum distance of a point $p_0$ to a submanifold $S$ of Euclidean space is quite a difficult problem in general. In the case that every intersection of $S$ with a ball is a closed set, we know that the distance function measuring the distance of points on $S$ to $p_0$ has a minimum because it is continuous. 
\skipall{
Necessary condition: Distance smaller than first focal distance.
Sufficient condition: Distance smaller than distance to cut locus
Distance of manifold to cut locus: Minimum of curvature radii and double normal half lengths
}
}

\section{Local Topological Structure of the SSB}
\label{secTopo}
\begin{figure}[tp]
\fbox{
	\parbox{7.5cm}{
		\center
		\includegraphics[width=7cm]{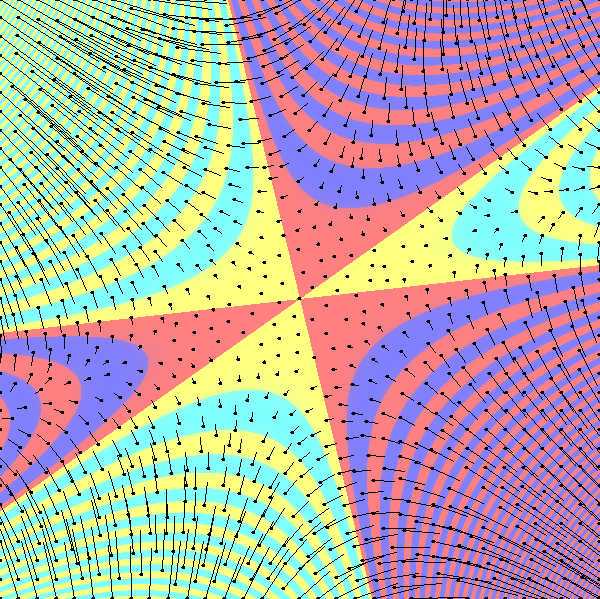}
		\caption{A homogeneous polynomial: Sign, isolines and gradients}
		\label{figSsbGrad} 
	}
}
\end{figure}

We want to provide some basic observations relevant in the context of discussing the local topological structure of the SSB, in particular in voltage space. It is helpful to know that if at some point $v_0$ in voltage space, the determinant of the differential $\D F(v_0)$ does not equal $0$, then the gradient of $\det \D F(v_0)$ is non-zero as well. Proof: Take the directional derivative of the $\D F$ at $v_0$ in the direction $\frac{v_0}{\abs {v_0}}$:
\begin{align}
\left. \dquot{\det \D F(t \cdot v_0)}{t} \right \rvert_{t=1} &= \left. \dquot{\det \begin{pmatrix}
\cdots & 2 t A_i v_0& \cdots 
\end{pmatrix}\transp
}{t}\right \rvert_{t=1} \nonumber \\
& = 
\left. \dquot{t^n \det \begin{pmatrix}
\cdots & 2 A_i v_0 &\cdots
\end{pmatrix}\transp}{t} \right \rvert _{t=1} \nonumber \\
& = 
\left. \dquot{t^n \det (\D F(v_0))}{t} \right \rvert _{t=1} \nonumber \\
& = 
\left. n t^{n-1} (\D F(v_0))\right \rvert _{t=1} \nonumber \\
& = n (\D F(v_0)) \nonumber \\
& \neq 0.
\end{align}
Of independent interest may be the observation that when moving along a ray from the origin uniformly parameterized by $t$, the determinant of the differential is scaled as $t^n$ compared to its value at $t=1$. 
This provides the intuition behind the proof: Since the function value scales with distance from the origin, unless the value is zero there will always be a radial component of the gradient. This is true for all homogeneous polynomials. \vref{figSsbGrad} shows the isolines of the bivariate cubic homogeneous polynomial $x^2y + y^2x + 0.1x^3-0.3y^3$ together with the direction and magnitude of its gradient field at selected points.

Another observation is that, in case that the rank of $\D F$ at some point is $n-2$ or lower, not only the determinant of $\D F$, but also its gradient vanish. Proof: $\grad \det \D F$ is a vector with components composed of partial derivatives in the basis directions $V_i$:
\begin{align}
\dquot{\det \D F}{V_i} & = \dquot{\det \begin{pmatrix}\dquot{F}{V_1} & \cdots & \dquot{F}{V_n} \end{pmatrix}}{V_i} \\
& = \sum_{l=1}^n\det \begin{pmatrix}\cdots & \dquot{F}{V_{l-1}} & \ddquot{F}{\der V_i\der V_l} & \cdots\end{pmatrix} .
\end{align}
All terms of the latter sum remain zero, since only one column of the original differential $\D F$ need to be changed to arrive at $\begin{pmatrix}\cdots & \dquot{F}{V_{l-1}} & \ddquot{F}{\der V_i\der V_l} & \cdots\end{pmatrix}$. Changing a single column can raise the rank of a matrix at most by one, and since we assumed the rank of $\D F$ was at most $n-2$, the changed matrices are still singular.

The converse also holds: If the rank of $\D F(v)$ is $n-1$, and thus $\det \D F(v)=0$, then the gradient of $\det \D F(v)$ is nonzero.
Proof: There is a direction $\delta$ so that $\det \D F(v+t\cdot \delta) \neq 0$ for small $t\in[0, \epsilon \in\IR^+]$. 
Thus the columns of $\D F(v+t\cdot \delta)$ span an $n$-dimensional volume $V(t) \eqdef \det \D F(v)$, which is nonzero if $t\neq 0$.
Without loss of generality, the last $n-1$ columns of $\D F(v+t\cdot \delta)$ can be assumed to span an $n-1$ dimensional volume $A(t)$ with $A(0)\neq 0$. Should this not be the case, we can always rearrange the columns to make it so. 
Let $h(t)$ be the distance between the hyperplane containing that volume and the point with coordinates given by the first column of $\D F(v+t\cdot \delta)$. Then we have $V(t) = h(t) \cdot A(t)$. Differentiating and setting $t=0$ we get 
\begin{align} \label{d0conv}
\dot V(0) = \dot h(0) \cdot A(0) + h(0) \cdot \dot A(0).
\end{align}
Because $0=V(0)=h(0)\cdot A(0)=0$ but $A(0)\neq 0$, it must be that $h(0)=0$. 
Hence we can simplify \vref{d0conv} to be $\dot V(0) = \dot h(0) \cdot A(0)$. That $\dot V(0)\neq 0$ implies that the gradient of $\det \D F(v)\neq 0$, because $\dot V(0)=\delta \cdot \grad \det \D F(v)$. 
Thus it remains to be shown that $\dot h(0)$ is nonzero. By rearranging the columns of $\D F(v+t\cdot \delta)$, we can always ensure that $\dot h(0)$, which is the speed with which the first column of deviates from the hyperplane spanned by the rest, is nonzero: 
There has to be at least one column that moves out of the plane spanned by all columns as $t$ changes, in order to ensure $V(t)\neq 0$ for $t \neq 0$. 
So we not only need to choose the first column so that the remaining $n-1$ columns have full rank, but also that it moves out of the $n-1$ dimensional subspace spanned by all columns fast enough (i.e. not quadratically or slower, but with a linear term). 
This always works because we assumed $V(t)\neq 0$ or $t\neq 0$, and in order to generate some nonzero volume $V(t)$, one point must move out of the plane and the others must span a nonzero ``area'' $A(t)$ for small enough $t$. This motion has a nonzero linear term because $D F(v)$ is linear in $v$, and therefore the derivative does not vanish.
%Note that to know $\dot h(0)$, the normal of that plane is needed, but not its derivative.
These considerations can also be summarized by observing that $\D F(v)$ being a linear function of $v$ and having a rank of $n-1$ means that the determinant of $\D F(v+t\cdot \delta)$, as a polynomial in $t$ of degree $n$, has a \emph{simple} root at $0$ for some choice of $\delta$, which implies that the derivative in the direction $\delta$ is nonzero, and therefore the gradient is nonzero.

\skipall{
1  0  0 
0 -3  0
0  0  1

0  0  3
0  2  0
3  0  0

2  0  0 
0 -5  0
0  0  2

@(1, 0, 1)

}
\subsection{Singularities on the SSB in Voltage Space}
\begin{figure}[tp]
\fbox{
	\parbox{7.5cm}{
		\center
		\includegraphics[width=7cm]{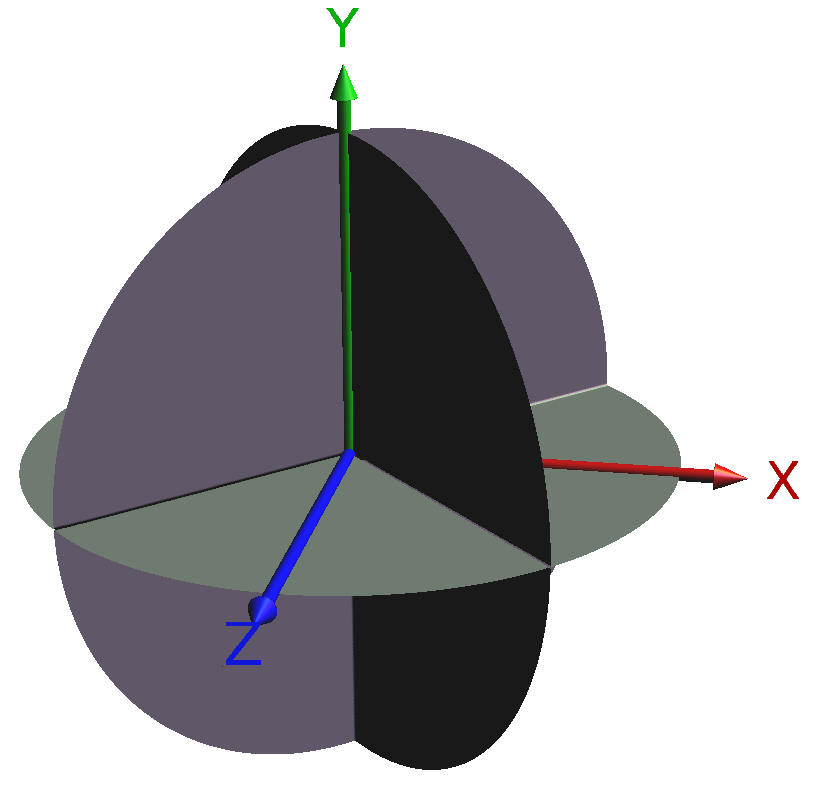}
		\caption{A determinant isosurface with non-manifold lines}
		\label{figDetDegen0} 
	}
}
\end{figure}

\begin{figure}[tp]
\fbox{
	\parbox{7.5cm}{
		\center
		\includegraphics[width=7cm]{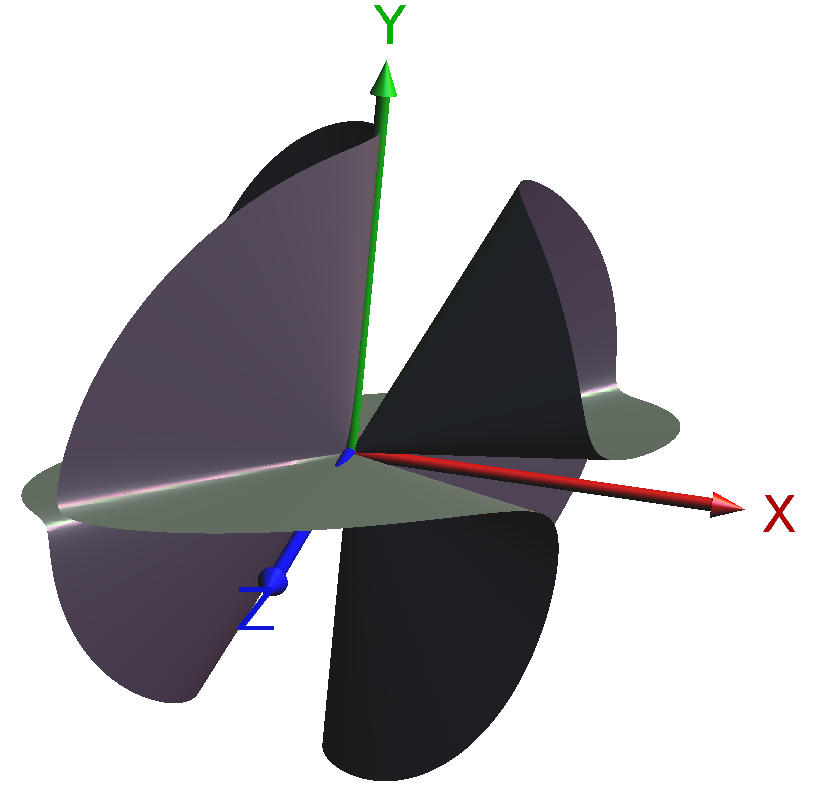}
		\caption{A perturbed determinant isosurface without non-manifold lines}
		\label{figDetNondegen0} 
	}
}
\end{figure}

Although it is not a generic case, it may happen that the $(\det \D F=0)$-surface for a quadratic function $F$ contains non-manifold points other than the origin. As an example consider the function
\begin{align}
F(v) = \begin{pmatrix}v\transp A_1 v \\ v\transp A_2 v \\ v\transp A_3 v \end{pmatrix}
\end{align}
where
\begin{align}
A_1 = \begin{pmatrix}
1&0&0\\0&-3&0\\0&0&1
\end{pmatrix},\,
A_2 = \begin{pmatrix}
0&0&3\\0&2&0\\3&0&0
\end{pmatrix},\,
A_3 = \begin{pmatrix}
2&0&0\\0&-5&0\\0&0&2
\end{pmatrix}.
\end{align}
Its differential is
\begin{align}
\D F\begin{pmatrix}x\\y\\z\end{pmatrix} = 2 \cdot 
\begin{pmatrix} 
x & -3y & z \\
3z & 2y & 3x \\
2x & -5y & 2z
\end{pmatrix},
\end{align}
which, when evaluated at $v\transp=(x,y,z)=(1,0,1)$, yields the rank-1 matrix
\begin{align}
2\cdot\begin{pmatrix} 
1 & 0 & 1 \\
3 & 0 & 3 \\
2 & 0 & 2
\end{pmatrix},
\end{align}
It is easily checked that at $(x,y,z)=(1,\epsilon,1)$ and $(x,y,z)=(1+\epsilon,0,1-\epsilon)$ for $\epsilon\neq0$, we generally obtain a rank-2 matrix.
Thus there are two directions (ignoring the trivial radial direction) in which the point $(x,y,z)=(1, 0, 1)$ can be perturbed so that $\D F$ stays degenerate, and thus directions that point along the SSB. However, changing $v$ by a linear combination of these directions generally results in a second-order deviation of $\det \D F$ from zero. Hence, the SSB is not locally a manifold at $v$. In fact, it looks as shown in \vref{figDetDegen0}: There are three radial lines composed of non-manifold points. This is bad for our intended use of differential geometric methods, because a tangential space is not defined at these points. This manifests itself in the computations by the gradient being zero, as explained in the previous subsection.
Luckily, this bad situation is not stable under perturbations of $F$. By a little change of two numbers in the definition of $A_1$, namely e.g. 
\begin{align}
A_1 \eqdef A_1 + \begin{pmatrix}
0&\epsilon&0\\
\epsilon&0&\frac\epsilon 2\\
0&\frac\epsilon 2&0
\end{pmatrix}, 
\end{align}
the SSB becomes a manifold (except at the origin), depicted in \vref{figDetNondegen0} for $\epsilon=0.1$. 

We prepared an animation of the evolution of the SSB as $\epsilon$ changes from $0$ to $0.6$: See the ancillary file {\tt{SSBTopologyChanges.mp4}} or \url{https://www.dropbox.com/s/9relj42dpu9i5re/SSBTopologyChanges.mp4}, first presented in the talks \cite{wolter2018NTU,wolter2018CGI}. This illustration reveals some other aspects of pathological SSBs:
\begin{itemize}
\item During the animation, two other degenerate situations occur. Both are associated with a change in topology: At the first such situation, a cone disappears. At the second, a fold in the single remaining sheet spawns off a cone. This illustrates that the set of all degenerate SSBs has codimension $1$ inside the set of all possible SSBs: On a curve through the space of all SSBs that connects two SSBs of different topology, there necessarily has to be a point where the SSB is non-manifold at points other than the origin.
\item The disappearing cone shows that the degeneracies where two portions of the SSB cross are not the only ones. To classify the possibilities in general, we need to look at the leading terms of the Taylor expansion around the degenerate point for the restriction of $\det \D F$ to a sphere around the origin. In our case, the dominant Taylor terms are all second order, forming a Hessian matrix. It is easy to see that an indefinite Hessian will result in intersecting SSB parts (Because when going around the singularity, the sign has to change four times) whereas a definite Hessian leads to 1-dimensional SSB parts that mark a cone's transition from being very narrow to nonexistent. In higher dimensions, the SSB can self-intersect in more complicated ways. In even less generic cases, the leading Taylor terms might then be higher order than quadratic.
\item The SSB at $\epsilon=0$, picured in \vref{figDetDegen0}, is the union of three planes. 
This is because the polynomial $\det \D F$ can be factorized into three linear terms, and the zero set of a product of functions is the union of their individual zero sets. One might be tempted to hope that the degenerate SSBs always can be decomposed as a union of simpler cone structures, and thus it may be possible to factorize the polynomial. But the self-intersecting SSB that occurs when the cone is spawned in the animation falsifies this, as here the surface clearly intersects with \emph{itself} by forming the loop that is about to become the new cone.
\end{itemize}

\begin{figure}[htp]
\fbox{
	\parbox{7.5cm}{
		\center
		\includegraphics[width=7.4cm]{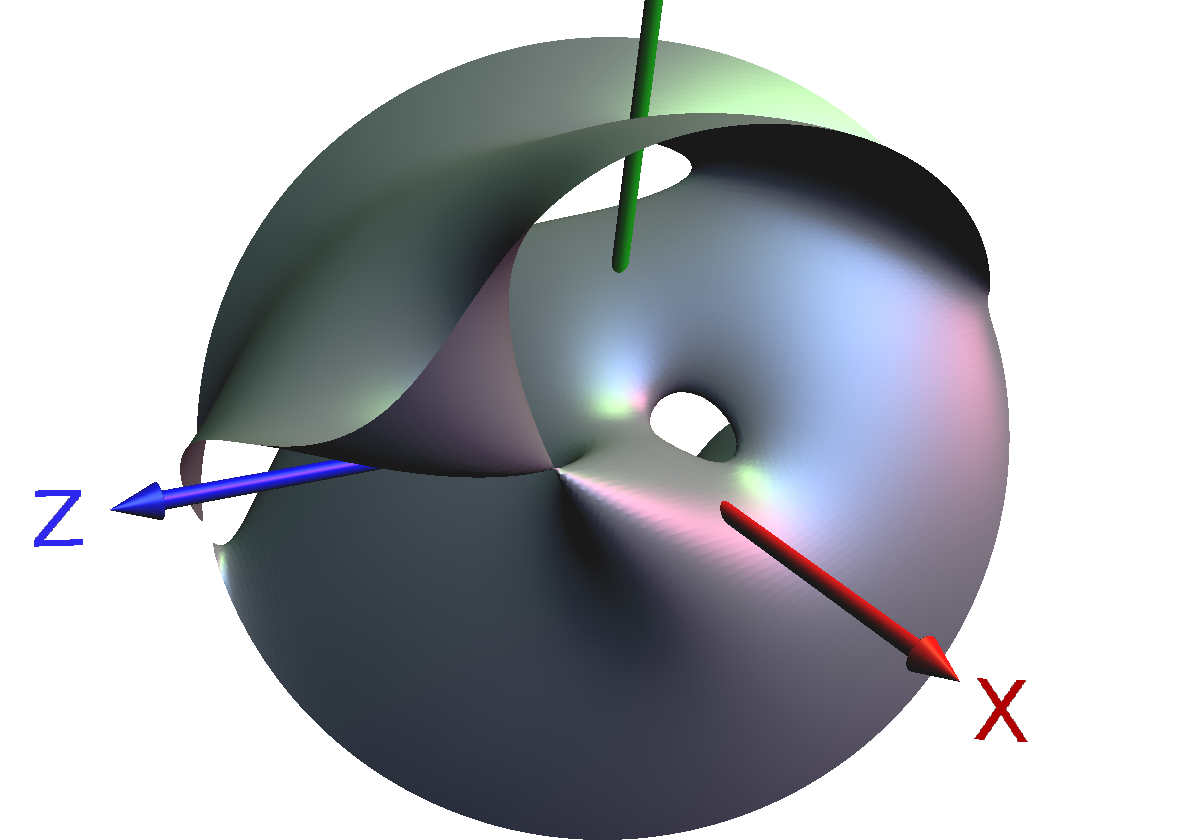}
		\caption{The intersection of an arbitrary SSB in 4 dimensions and half of the unit 3-sphere, stereographically projected into 3-space.}
		\label{figSSB4DevoFrame} 
	}
}
\end{figure}

Note that these depictions of the SSB are highly redundant: Since it is composed of radial lines and has point symmetry about the origin, an SSB for $n=3$ can be described completely by its intersection with (half) a unit $2$-sphere. The degenerate situations then manifest on this sphere as isolated points and crossing lines. One might think that the generalization to higher $n$ (say, $n=4$) then involves isolated lines and crossing planes on the $(n-1)$-sphere. This however is not the generic situation for locations where the rank of the differential is less than $n-1$. Rather, the usual case looks like the zero set of a quadratic form, which for $n=4$ means isolated points and double cones. To give an idea about the possible SSB topologies for $n=4$ and their changes as the coefficients of the quadratic forms that define the power flow map move along a line in configuration space, we have prepared two animations, available under \url{https://www.dropbox.com/s/go29bk47w5nrazw/SSB4D_1.mp4} and \cite{wolter2018NTU,wolter2018CGI}.
\url{https://www.dropbox.com/s/jo7xoo07riicd7w/SSB4D_2.mp4} or as ancillary files of this paper, named {\tt{SSB4D_1.mp4}} and  {\tt{SSB4D_2.mp4}} and first presented at \cite{wolter2018CGI}.
The videos show the stereographic projection of those parts of the evolving SSBs that intersect with one half of the 3-sphere. The coefficients and the directions for changing them have been chosen randomly. 
\vref{figSSB4DevoFrame} shows a frame from the video exhibiting a typical double-cone degeneracy.

\begin{figure}[htp]
\fbox{
	\parbox{7.5cm}{
		\center
		\includegraphics[width=7.4cm]{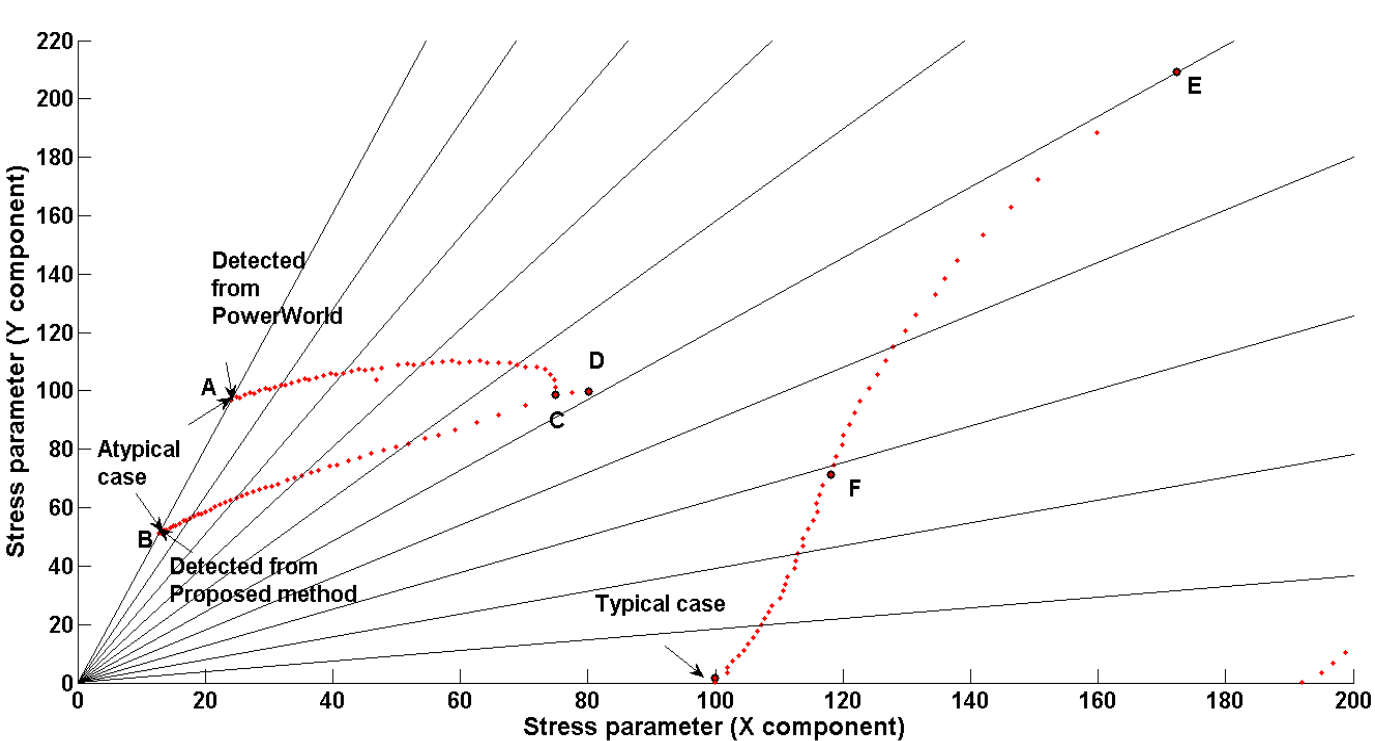}
		\caption{Figure 26 in \cite{makarov2014non}. Reproduction of the graphic should be fair use both because it is for purposes of criticism and because this is scholary work, see 17 U.S. Code § 107.}
		\label{figWrong} 
	}
}
\end{figure}

We were motivated to investigate the possible local topology of SSBs not only because our differential geometric algorithms require it to be smooth, but also because Figure 26 in \cite{makarov2014non}, reproduced here in \vref{figWrong}, depicts a $2$-dimensional section through an SSB that seems to exhibit 3-way branching (Point C in that figure) and a part of the SSB simply ending (Point D). We think that this should not be possible, and will now prove that any point in a 2D-section of an SSB should be surrounded by an even number of outgoing lines. The features of point C with three outgoing lines and point D with a single line can possibly be explained by either the algorithm of \cite{makarov2014non} detecting SSB points between C and D where there are none and instead $\det \D F$ is merely very small, or (less likely, given how their algorithm works) the algorithm missing a fourth branch coming out of at C downwards and also not seeing the continuation of the line at D, which could be due to the appearance of another line at point E further back on the same stress direction ray.

The claim that every point in a 2D-section of the zero set of an algebraic curve should be surrounded by an even number of outgoing lines has already been made in some form in the dissertation of Gauß. We refer to \cite{ghys2016singular} for historical context and a proof of the claim.
The intersection of the SSB with any affine 2D subspace of $\IR^n$ is an algebraic curve as the SSB can be defined as the zero set of the multivariate polynomial $\det \D F$, and the restriction of a multivariate polynomial to an affine subspace with a basis is again a polynomial in the variables given by the basis directions of the subspace. Hence, the theorem proven in \cite{ghys2016singular} is applicable to the case at hand.

%\section{Using geodesics for computing paths on implicit submanifolds}

\section{Improving Optima Under Security Constraints}
Optimal power flow (OPF) is a subfield of power flow studies that deals with the problem of finding values of controllable variables such as power generation and generator voltage magnitudes that are optimal with respect to some goal function such as generation cost or network losses. The resulting optimal operating point should also be constrained not to lie too close to the SSB (in power space) so as to minimize risk of voltage instability. Furthermore, one may want to consider contingencies where buses or transmission lines fail, leading to altered network topology and therefore a different SSB. So the operating point should also keep a safe distance to all considered contingency SSBs. Additionally, the variables may be subject to inequality constraints corresponding to engineering limits of the hardware.

Here we will present an algorithm for ensuring a minimum distance to the SSB through power space, given an initial optimal operating point. This is accomplished by ``pushing away'' the operating point from the SSB while maintaining local optimality. The procedure can easily be adapted to push the point away from multiple SSBs simultaneously. We do not address the question of how to find the initial operating point, only how to improve its security.

Similar to the above algorithm for finding perpendicular foot points, our algorithm is a combination of a predictor step given by an ODE, which may occasionally fail, and a corrector step which usually has nothing to do except if inaccuracies accumulate or the predictor fails. The corrector step essentially re-solves the optimization problem, given the output of the predictor as the starting point. It is also possible to omit the predictor altogether and just feed the last value of the operating point to the corrector as the initial value for computing the next time step.

\begin{figure*}[!htp] 
{
	\parbox{\textwidth}{
		\center
		\subfigure[The optimal operating point, subject to 2 linear inequality constraints, and its closest point on the boundary. Initially, only one constraint is active.]{\includegraphics[width=0.3\textwidth]{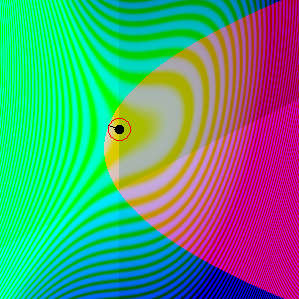}}\quad
		\subfigure[The optimum is pushed away from the boundary while the first constraint remains active.]{\includegraphics[width=0.3\textwidth]{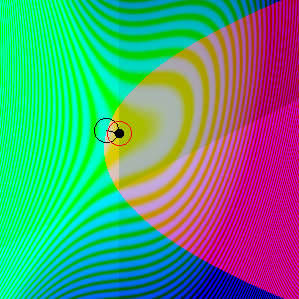}}\quad
		\subfigure[When the optimum is pushed further, both constraints are inactive.]{\includegraphics[width=0.3\textwidth]{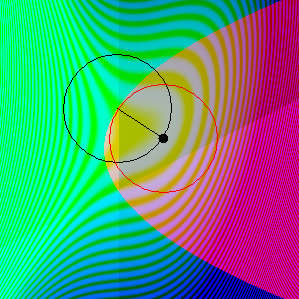}}\\
		\subfigure[A second closest point on the boundary is found. The optimum now moves along the cut locus.]{\includegraphics[width=0.3\textwidth]{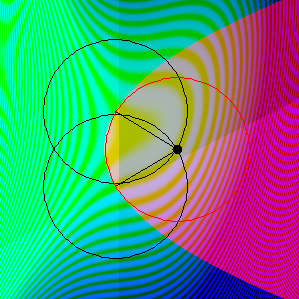}}\quad
		\subfigure[The second constraint becomes active, causing one of the minimum-distance-to-boundary constraints to become inactive.]{\includegraphics[width=0.3\textwidth]{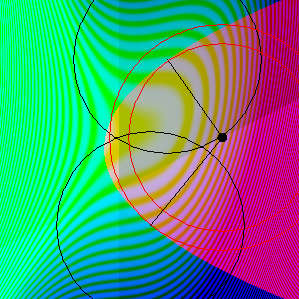}}\\
		\caption{Low-dimensional example run of the algorithm pushing the optimal operating point away from the boundary. Minimum distance constraints are shown using black circles and other constraints are visualized by using darker colors in the permissible region. Isolines of the goal function are also shown.}
		\label{figOPFA} 
	}
}
\end{figure*}

In order to push the operating point (OP) away from the SSB, it suffices to know all SSB points that are closest to the OP and push the OP in a direction that increases the distance from these, while accounting for the shifting of these closest point as the OP moves. When this pushes the OP across the cut locus of the SSB, it becomes necessary to add another closest point to consideration. Thus the corrector step includes a search for closest points and may change the dimension of the ODE when it discovers a new one. \vref{figOPFA} illustrates the idea, but note that the equations underlying the particular problem depicted have nothing to do with power flow per se.
\subsection{Variables and Constraints}

The time variable of the ODE will be called $t$. All other variables depend implicitly on it. At time $t$, the algorithm ensures that the OP is at least $t-\epsilon$ units away from the SSB, where $\epsilon$ is the maximum step size.

The operating point is given by the adjustable variables that are optimized for (such as power injections and squared voltage magnitudes), as well as the voltage phasors. The non-adjustable power variables are assumed fixed and folded as constants into the formulas. Let us call the variables $o_i$ for the adjustable variables and $v_j$ for the voltage variables, for suitable ranges of $i$ and $j$. Let $M(o, v)$ be the the difference vector of the squared voltage magnitude variables and power variables on the one side, and the squared voltage magnitudes and power injections as computed from the voltage variables $v$ according to the power flow equation on the other side. Then we have equality constraints $M(o, v)=0$ that govern the relation between the $o$-variables and the $v$-variables. If there are multiple SSBs to consider, we will need one set of $v$-variables for each one, along with a corresponding set of constraints.

We will continue to denote the power space coordinates of the OP with $p_i$ for $1\leq i\leq n$. Some of these may be constant and others may be just other names for an $o$-variable.

For each inequality constraint placed on the operating point, we need a slack variable $s_i$. From our understanding of the OPF problem, inequality constraints will not depend directly on voltage variables, so each inequality constraint can be represented by an equality constraint $L_i(o)- {s_i^2} = 0$, where $L_i$ is a differentiable function that has to be positive for the constraint to be satisfied.

The goal function that is to be minimized is called $G^{\text o}$. Its arguments are the $o$ variables.
The Lagrangian for the main optimization system is called $\mathcal L^{\text o}$. Its parameters are the $o$ and $v$ and slack variables and some Lagrange multipliers. Some of the variables mentioned in the following enter into it, but are not part of the variables that the main optimization problem optimizes over.

For each closest point we will need a set of variables describing that point. The closest points will be called ``feet'' from now on because they are the foot points of perpendiculars from the OP onto the SSB. The algorithm maintains a set of feet, some of which may be globally closest to the OP while others may be farther away. Should a foot be too close to another foot, or too far away from the OP, it should be dropped from the foot set to avoid needless computations. The foot with index $i$ is characterized by
\begin{itemize}
\item A point in voltage space, given by $n$ voltage variables $w_{ij}$ for $1\leq j\leq n$.
\item The kernel of $(\D F)\transp$ at $w_i$, given by $n$ variables $k_{ij}$ for $1\leq j\leq n$. The kernel is also the normal on the SSB in power space.
\item A constraint $\abs{k_i}=1$ to ensure that the kernel is nonzero, and constraints $(\D F(w_i))\transp k_i=0$ that ensure that $k_i$ really is the kernel. These constraints make sure that $w_i$ is a point on the SSB.
The first constraint can also be written as $\abs{k_i}^2=1$, which is not as good for the numerics of the corrector algorithm but simplifies the algebraic manipulations needed for the predictor. We can use one form of the constraint in the predictor and the other form in the corrector without problem.
\item A slack variable $\sigma_i$.
\item A constraint $\abs{k_i \cdot (p - F(w_i))} - t -  {\sigma_i^2} = 0$ that ensures that the distance between the foot and the OP is at least $t$.
\skipall{\item Constraints of the form $(p - F(w_i)) = \shin_i k_i$ (with an extra variable $\shin_i$) or similar can also be used to provide more stability and faster convergence by explicitly expressing that the direction from the foot to to OP shold be perpendicular to the SSB. We will not use these in the following for brevity. }
\item A goal function $G^{\text f}_i(o, w_i) = \abs{F(w_i) - p}^2$ that is used to minimize the distance between the OP and the foot.
\item A Lagrangian $\mathcal L^{\text f}_i$ for the optimization problem just mentioned. Its parameters are the $w_i$, $k_i$ and $\sigma_i$ variables and some Lagrange multipliers. Some of the $o$ variables occur in $\mathcal L^f_i$ (via the $p$ in the goal function), but are not part of the variables that the optimization problem for this foot optimizes over.
\end{itemize}

For each constraint we will need a Lagrange multiplier. We will call them:
\begin{itemize}
\item $\lambda^{\text M}_i$ for the $i$-th line of the constraint equation system $M(o, v) = 0$
\item $\lambda^{\text L}_i$ for the constraint $L_i(o)- {s_i^2} = 0$
\item $\lambda^{\text N}_i$ for the constraint $\abs{k_i}=1$ or $\abs{k_i}^2=1$
\item $\lambda^{\text k}_{ij}$ for the $j$-th line of the constraint equation system $(\D F(w_i))\transp k_i=0$
\item $\lambda^{\text d}_{i}$ for the constraint $\abs{k_i \cdot (p - F(w_i))} - t - {\sigma_i^2} = 0$
\end{itemize}

\subsection{Predictor ODE}
Since the variables are by assumption fulfilling all constraints, they form a stationary point of the Lagrangian of the optimization problem. There are actually several optimization problems: One for the main goal function and one for each foot.
The Lagrangians are:
\begin{align}
\mathcal L^{\text o}(o, v, s, \sigma, \lambda^{\text M}, \lambda^{\text L}, \lambda^{\text d}) 
  & = G^{\text o}(o) + \sum_i \lambda^{\text M} \cdot M(o, v) + \nonumber \\
  &  \sum_i \lambda^{\text L}_i (L_i(o)- {s_i^2}) + \\
  & \sum_i \lambda^{\text d}_{i} (\abs{k_i \cdot (p - F(w_i))} - t -  {\sigma_i^2}) \nonumber\\
\mathcal L^{\text f}_i(w_i, k_i, \lambda^{\text N}_i, \lambda^{\text k}_{i}) & = \abs{F(w_i) - p}^2 
   + \lambda^{\text N}_i (\abs{k_i}^2-1) + \nonumber\\
 &   k_i\transp(\D F(w_i)) \lambda^{\text k}_{i} .
\end{align}
The conditions that these Lagrangians be stationary with respect to their explicit parameters are then
\begin{align}
0 &= \grad G^{\text o}(o) + \lambda^{\text M}_j \grad_o M(o, v)_j \\
  &+ \sum_i \lambda^{\text L}_i (\grad L_i(o))_j \nonumber \\
  &+ \sum_i \lambda^{\text d}_{i} \sign({k_i \cdot (p - F(w_i))}) k_{ij} \nonumber\\
0 &= \sum_i \lambda^{\text M}_i \grad_v M(o, v)_i \\
0 &= \lambda^{\text L}_i {s_i} \\
0 &= \lambda^{\text d}_i {\sigma_i} \\
0 &= M(o, v) \\
0 &= L_i(o)-{s_i^2}\\
0 &= \abs{k_i \cdot (p - F(w_i))} - t - {\sigma_i^2}\\
0 &= 2(F(w_i)_j - p_j) + k_i\transp \left (\dquot{\D F(w_i)}{w_{ij}}\right ) \lambda^{\text k}_{i} \\
0 &= \lambda^{\text N}_i {k_i} \\
0 &= \abs{k_i}^2-1 \\
0 &= (\D F(w_i))\transp k_i .
\end{align}
Here, an equation containing an unbound index variable is actually several equations, one for each possible value of the index variable.
Note that we ignore the null set where the derivative of the absolute and sign functions at $0$ would be needed. Any errors that might come from that will be mended by the corrector step.

Next, we differentiate all these equations with respect to $t$ to express that the variables remain a solution of the optimization problems as $t$ changes. Again, we express differentiation with respect to $t$ with a dot on top of symbols. The equations are then:
\begin{align}
0 &= (\Hs G^{\text o}(o)) \dot o \\
  &+ \dot \lambda^{\text M}_j \grad_o M(o, v)_j\nonumber \\
  &+ \sum_i \lambda^{\text M}_i ((\Hs_o M(o, v)_i) \dot o)_j \nonumber\\
  &+ \sum_i \lambda^{\text M}_i (\grad_o ((\D_v M(o, v)_i) \dot v))_j \nonumber\\
  &+ \sum_i \dot \lambda^{\text L}_i (\grad L_i(o))_j \nonumber \\
  &+ \sum_i \lambda^{\text L}_i ((\Hs L_i(o)) \dot o)_j \nonumber \\
  &+ \sum_i \dot \lambda^{\text d}_{i} \sign({k_i \cdot (p - F(w_i))}) k_{ij} \nonumber\\
  &+ \sum_i      \lambda^{\text d}_{i} \sign({k_i \cdot (p - F(w_i))}) \dot k_{ij} \nonumber\\
0 &= \sum_i \dot \lambda^{\text M}_i (\grad_v M(o, v)_i)_j \\
  &+ \sum_i \lambda^{\text M}_i (\grad_v ((\D_o M(o, v)_i) \dot o))_j \nonumber\\
  &+ \sum_i \lambda^{\text M}_i ((\Hs_v M(o, v)_i) \dot v)_j \nonumber\\
0 &= (\dot \lambda^{\text L}_i {s_i} + \lambda^{\text L}_i {\dot s_i})\\
0 &= (\dot \lambda^{\text d}_i {\sigma_i} + \lambda^{\text d}_i {\dot \sigma_i}))\\
0 &= (\D_o M(o, v)) \dot o + (\D_v M(o, v)) \dot v\\
0 &= (\D L_i(o))\dot o - 2 {s_i} \dot s_i\\
1 &= \sign({k_i \cdot (p - F(w_i))}) \\
  & \quad \cdot (\dot k_i \cdot (p - F(w_i)) + k_i \cdot (\dot p - (\D F(w_i))\dot w_i) ) \nonumber\\
  & - 2{\sigma_i}  \dot \sigma_i \nonumber\\
0 &= 2((\D F(w_i))_j \dot w_i - \dot p_j) \\
  &+ \dot k_i\transp \left (\dquot{\D F(w_i)}{w_{ij}}\right ) \lambda^{\text k}_{i}\nonumber \\
  &+ k_i\transp \left (\dquot{\D F(w_i)}{w_{ij}}\right ) \dot \lambda^{\text k}_{i}\nonumber \\
\intertext{(Note that in the previous equation, a term containing third derivatives has been dropped as $F$ is quadratic. Apart from that, these equations should work in similar applications where $F$ is not quadratic.)}
0 &= \dot \lambda^{\text N}_i {k_i} + \lambda^{\text N}_i {\dot k_i} \\
0 &= k_i \cdot \dot k_i \\
0 &= ((\Hs F(w_i))\dot w_i)\transp k_i + (\D F(w_i))\transp \dot k_i .
\end{align}
This is a linear equation system for the unknowns $\dot o$, $\dot v$, $\dot s$, $\dot \sigma$, $\dot w$, $\dot k$, 
$\dot \lambda^{\text M}$, $\dot \lambda^{\text L}$, $\dot \lambda^{\text N}$, $\dot \lambda^{\text k}$ and $\dot \lambda^{\text d}$.
Its solution can be used in a numerical ODE integrator to predict how the the OP and the foot as well as the slack variables and Lagrange multipliers change while $t$, the minimum distance to the SSB, increases to whatever value desired (and feasible). If the predicted stationary point ceases to be a local optimum, the corrector will step in and fix the situation by choosing a nearby local optimum, if one exists.
Note that in typical applications, $\Hs G^{\text o}$ and $(\D L_i)$ will be sparse, so the whole equation system is sparse.
\subsection{Corrector Step}
In the corrector step, the optimization problems are re-solved for the current value of $t$ using the variable values obtained from the predictor as the initialization. Also, an attempt is made to discover new feet. 

We use a kind of projected gradient descent for solving the optimization problems. 
All constraints can be written in the Form $c_i(x)=0$, where $x$ are the variables that are being optimized and $c_i$ is a function associated with the $i$-th constraint. Let $\tilde c_i$ be the gradient of $c_i$, at least for now. Then in each step the solver determines the space spanned by all the $\tilde c_i(x)$, where $x$ is the current state of the variables. Within that space, it determines for each $i$ the hyperplane formed by those $y$ with $c_i(x) + y\cdot \tilde c_i(x) = 0$, and then finds the point $z$ with the least squared distance to these hyperplanes. This should converge to the set of simultaneous solutions for all constraints similar to Newton iterations; indeed, Newton iterations are a special case of this. Let $\tilde G$ be the gradient of the goal function $G$, at least for now. Let $\delta$ be the orthogonal projection of $\tilde G(x)$ on the orthogonal complement of the space spanned by the $\tilde c_i(x)$. Then $x$ is updated in each step according to $x \eqdef z + \gamma \delta$, where $\gamma$ is an adaptively chosen step size. This is repeated until $\abs\delta$ and all $c_i(x)$ are small enough or the maximum number of iterations is reached.

Note that the constraints contain squared slack variables $s_i^2$ and $\sigma_i^2$. When these are small or zero, but the solution requires them to have a nonzero value, there will be a problem because the corresponding gradient component of the associated $c_i$-function is small or zero too. So for the corrector, it is better to switch to absolute slack variables $\abs{s_i}$ and $\abs{\sigma_i}$. To decide on a direction of the gradient when the slack variable is zero, we use the convention $\left. \dquot{\abs{a}}{a}\right\vert_{a=0} = 1$.

We have several interdependent optimization problems to solve. The most rigorous approach would be to re-solve the optimization problems that adjust the feet so as to minimize the foot-OP distances in each step of the the main optimization problem. This means that inside the main optimizer loop, there is a nested loop for the foot optimization. However, we devised a way to use a single optimizer loop. The idea is that minimizing the distance of a foot to the OP should not shift the OP, and enforcing the required distance constraint between a foot and the OP should not move the foot. What this means in practice is that if $G=G^{\text f}_i$, we do not use the gradient for $\tilde G(x)$ but the gradient with the components belonging to OP variables set to zero. Likewise, for those constraints $c_i$ that ensure the minimum distance between the OP and a foot, $c_i(\dots, p, k_i, w_i, \sigma_i, \dots) = \abs{k_i \cdot (p - F(w_i))} - t - {\sigma_i^2} = 0$, we set to zero those components of $\tilde c_i$ that correspond to the foot-associated variables $ k_i$ and $w_i$. We call this scheme ``reactionless gradient descent'' because there are interactions between parts of the solution where an influencing ``force'' only acts on one part.

Reactionless gradient descent has the potential to be faster than the nested optimization loops because it solves several interrelated optimizations problems at once, but because the vector field that it follows to optimize the functions is not really a gradient field, convergence is not assured unconditionally. Also, in some cases it may converge, but rather slowly. We conjecture that in our case, the convergence properties have to do with the curvature of the SSB. If we investigate this further, we may be able to find a criterion to switch between conventional projected gradient descent with nested loops and reactionless gradient descent in order to get the best of both worlds.

After the corrector step, the new values for the Lagrange multipliers need to be determined in case the corrector is used together with the predictor, as the corrector does not involve Lagrange multipliers.

\subsection{Finding New Foot Points}
When the OP is pushed across the cut locus of the SSB, a new closest point on the SSB may arise. We aim to find this point soon enough for the algorithm to be correct. Solving the global minimum distance problem is hard, but it seems if we are only interested in a point that is within $\epsilon$ of the globally shortest distance, there is an easier way. We proceed by the following steps: 
\begin{enumerate}
\item Solve the nonlinear equation system with the two equations $\det \D F(y)=0$ and $\abs{F(y)-p} = \abs{F(z)-p} - \epsilon$ for $y$, where $p$ is the position of the OP in power space and $z$ is the closest currently known foot point. If no solution exists, continue at step 5. The first equation may be replaced with a suitable equivalent equation (system), such as requiring $\D F(y)$ to have a zero eigenvalue and associated kernel vector.
This has the advantage of making the equations quadratic and the next step trivial.
\item Find the left unit eigenvector $N_P$ of $\D F(y)$ for the eigenvalue $0$. 
\item Use projected gradient descent to update $y$ and $N_P$ so that $\abs{F(y)-p}$ is minimized under the constraints $N_P \cdot N_P = 1$ and $(\D F(y))\transp N_P = 0$. The constraints ensure that $y$ remains a point on the SSB and that $N_P$ remains the left eigenvector for eigenvalue $0$. The goal function leads to $y$ being a locally closest point.
\item Let $m$ be the first unused foot index. Add a new foot with index $m$, initialized with $w_m = y$, $k_m=N_P$ and the appropriate value for the slack variable $\sigma_m$. Go back to step 1.
\item If at least one foot has been added in the loop, make sure that the corrector is run afterwards.
\item Remove feet that are (near) duplicates of other feet with smaller index, or that are too far away from the OP compared to the closest feet so they likely won't be relevant in the near future.
\end{enumerate}
It is crucial for this algorithm that the equation system in step 1 can be solved reliably and efficiently. In our implementation, we use the same algorithm as for the projected gradient descent, albeit with a constant goal function so that it is effectively only executing the projection step. In low dimensional examples this works fine, but we have not yet tested if this holds up in higher dimensions.

\section{Conclusion and Further Work}

We presented various formulas and algorithms for calculating geometric entities associated with the SSB, such as curvatures, geodesics and orthogonal projections of points and curves. We laid special focus on the SSB in voltage space, which historically has received less attention and which is less accessible computationally, needing in general one order of derivatives more. 
Our algorithm for local inversion of the power flow map allows us to continue a solution to the power flow equations along a curve with high precision. 

We also clarified some aspects of the local topological structure of the SSB. Furthermore, we presented an algorithm for improving the results of an optimal power flow computations so as to fulfill certain security constraints a posteriori. The reactionless gradient descent algorithm warrants further study; here, understanding the curvature of the SSB in power space might play a role in its analysis.

We showed how to compute the exponential map and its differential on the SSB in voltage space. This is a preparation for future work wherein we want to exploit this for finding points (with some coordinates specified) on a submanifold of the SSB using an algorithm from \cite{gutschke2015differential}. This in turn should be useful in finding safe paths through state space for restoring the power system after a contingency.

Hoping that we could motivate the use of differential geometric methods for power flow related problems, we look forward to providing more solutions employing applications of these methods.

\subsection{Acknowledgements}
This work had its origin during a research term of F.-E. Wolter at MIT in summer 2016, supported by a MISTI (MIT Germany) research grant that he initiated together with T. Sapsis, N. Patrikalakis and S. Karaman. This group was later joined by K. Turitsin. During this research term, F.-E. Wolter came into contact with electrical power researchers from Argonne Laboratory (D. Molzahn) and the University of Michigan (Ian Hiskens), who later invited him to visit. 

All these people have been supporting the advancement of the ideas presented in this paper via various discussions and references, which had been both enlightening and encouraging. Later on during visits to the IMI research institute in Singapore which included presenting seminars there \cite{wolter2018NTU,wolter2018CGI}, work done in this project evolved further. Most recently, F.-E. Wolter, during his 2019 visit in Singapore \cite{wolter2019NTU}, had the chance to come into contact with power flow researcher H. D. Nguyen who provided stimulating questions. 

This is a good occasion to say thank you to the IMI for the opportunity to present ideas of this project in interdisciplinary seminars and for many stimulating informal discussions with various researchers there, as well as of course for the the generous invitations by N. Thalmann that helped building new promising cooperations.

\bibliographystyle{elsarticle-num}
\bibliography{Paper}
%\section{Appendix}
\end{document}